\documentclass{amsart}
\usepackage{amsmath,amsthm,amsfonts,amssymb,epsfig,graphicx,verbatim}
\usepackage[usenames,dvips]{color}
\input xy
\xyoption{all}

\newtheorem{thm}{Theorem}[section] 
\newtheorem{lemma}[thm]{Lemma}

\newtheorem{proposition}[thm]{Proposition}

\theoremstyle{definition} 
\newtheorem{definition}[thm]{Definition}
\newtheorem{problem}[thm]{Problem}
\newtheorem{question}[thm]{Question}

\newtheorem{example}[thm]{Example}

\newtheorem*{ack*}{Acknowledgements}
\theoremstyle{remark} 
\newtheorem{remark}{Remark}

\newtheorem*{thm*}{Theorem}
\newtheorem*{lemma*}{Lemma}
\newtheorem*{tlm*}{Technical Lemma}
\newtheorem*{sblm*}{Sublemma}
\newtheorem*{corollary*}{Corollary}
\newtheorem*{proposition*}{Proposition}
\newtheorem*{criterion*}{Criterion}
\newtheorem*{algorithm*}{Algorithm}
\newtheorem*{conjecture*}{Conjecture}
\newtheorem*{conclusion*}{Conclusion}
\newtheorem*{answer*}{\textsf{Answer}}
\newtheorem*{example*}{Example}
\newtheorem*{remark*}{Remark}

\newcommand{\B}{\ensuremath{\mathcal{B}}}

\newcommand{\C}{\ensuremath{\mathbb{C}}}

\newcommand{\Z}{\ensuremath{\mathbb{Z}}}
\newcommand{\N}{\ensuremath{\mathbb{N}}}

\def\II{\mbox{\rm \Large 1\hskip -0.353em 1}}
\newcommand{\ds}{\displaystyle}
\newcommand{\fr}{\displaystyle\frac}
\newcommand{\ben}{\begin{enumerate}}
\newcommand{\een}{\end{enumerate}}

\newcommand{\bq}{\begin{question}}
\newcommand{\eq}{\end{question}}
\newcommand{\bp}{\begin{problem}}
\newcommand{\ep}{\end{problem}}
\newcommand{\ba}{\begin{answer*}}
\newcommand{\ea}{\end{answer*}}

\newcommand{\ess}{\overset{\text{ess}}{=}}

\newcommand{\freq}{\operatorname{freq}}

\begin{document}

\numberwithin{equation}{section}

\title{Dynamical properties of the Pascal adic transformation}

\date{\today}
\author{Xavier M\'ela}
\address{IML, 163 avenue de Luminy, Case 907, 13288 Marseille Cedex 09, FRANCE}
\email{mela@iml.univ-mrs.fr}

\author{Karl Petersen}
\address{Department of Mathematics, CB 3250 Phillips Hall, University of North Carolina, Chapel Hill, NC 27599, USA}
\email{petersen@math.unc.edu}


\begin{abstract}
We study the dynamics of a transformation that acts on infinite paths
in the graph associated with Pascal's triangle. For each ergodic
invariant measure the asymptotic law of the return time to cylinders
is given by a step function. We construct a representation of the
system by a subshift on a two-symbol alphabet and then prove that the
complexity function of this subshift is asymptotic to a cubic, the
frequencies of occurrence of blocks behave in a regular manner, and
the subshift is topologically weak mixing.
\end{abstract}

\maketitle


\section{Introduction}\label{sec:introduction}

Adic transformations were introduced by A. Vershik as combinatorial
models of the cutting and stacking constructions familiar in ergodic
theory  \cite{vershik7,vershik9,vershik5}. They move sequences
in a transverse manner to the usual shift transformation, much as the
horocycle flow is transverse to the geodesic flow
(cf. \cite{Sinai,Ito}). An adic transformation acts on the space of
infinite paths on an infinite graded graph, or Bratteli diagram, and
the dynamics of these transformations can provide information about
the associated $C^*$ algebras or families of group representations
(see for example
\cite{vershik1,vershik2,vershik3,vershik4}). Stationary adic
transformations (in which after the first, or root, level all levels
of the graph have the same number of vertices and the same pattern of
connections to the next level) correspond to odometers and
substitution subshifts
\cite{forrest, livshitz, vershik5, host1, solomyak2, solomyak3,durand-host-skau}.
 Every minimal homeomorphism of the Cantor set is
topologically conjugate to a particular type of adic transformation \cite{herman-putnam-skau},
and every ergodic measure-preserving transformation on a Lebesgue
space is isomorphic to an adic transformation with a unique nonatomic invariant
measure \cite{vershik9,vershik8}. The families of invariant sets for adic
transformations correspond to tail fields in probability theory and
ergodic theory, so ergodicity of invariant and quasi-invariant
measures for adic systems is equivalent to 0,1 laws, which guarantee
the triviality of these tail fields---see
\cite{gibbs,schmidt97,schmidt99,santiago}.
An especially regular and simple-looking nonstationary adic is the one
based on the Pascal triangle regarded as a graded graph.
Its $\sigma$-algebra of invariant sets corresponds to the exchangeable
or symmetric $\sigma$-algebra in $\{ 0,1\}^\mathbb N$, the sets fixed
by any permutation of finitely many coordinates (whereas the usual
tail $\sigma$-algebra consists of the sets invariant under any change
of finitely many coordinates).

In this paper we establish several dynamical properties of the Pascal
adic transformation. It is known that the set of nonatomic ergodic invariant measures
for this system is a one-parameter family
corresponding to the Bernoulli measures on $\{ 0,1\}^\mathbb N$, as explained below.
For each ergodic invariant measure we identify the asymptotic law of
the return time to cylinder sets determined by finite initial path
segments (Theorem \ref{thm:returntimes}). The original Pascal
transformation is not defined everywhere, which means that we are
dealing with a noncompact topological dynamical system. Attempts at
compactification or at forming quotients lead to discontinuities. To
overcome these difficulties, we use a countable family of
substitutions to produce a subshift on an alphabet of two symbols, $\{
a,b\}$, which represents the Pascal adic except for countably many
points (Theorem \ref{thm:pascal-sub-iso}).
This subshift consists of all subwords of all ``basic words'' formed
by concatenating words, rather than adding integers, in Pascal's
triangle. The basic word found at place $k$ in row $n$ has length
equal to the binomial coefficient $C(n,k)$ found at the same place in
the actual Pascal triangle, and the structure of the word conveys some
information about the history of its formation and therefore also some
extra information about the binomial coefficient.
Not only does this subshift have zero entropy (it was known before
that all the
invariant measures for the Pascal adic transformation have entropy
zero), but we can determine its complexity function: the number of $n$-blocks
is asymptotic to $n^3/6$ (Theorem \ref{thm:complexity}).
While the subshift supports uncountably many ergodic invariant
measures, it has a property that we call {\em directional unique
ergodicity}: once a ray in Pascal's triangle beginning at the root is
fixed, when we consider occurrences of a given block $B$ only in
basic blocks near that ray, the limiting frequency of occurrences
exists and equals the measure of the cylinder set $[B]$ according to
the ergodic invariant measure parametrized by the angle of the ray
(Theorem \ref{thm:freq}).
Finally, we use a characterization of weak mixing by Keynes and
Robertson \cite{keynes-robertson} and Weyl's theorem on uniform distribution to show
that the subshift is topologically weakly mixing (Theorem \ref{thm:pascal-twm}).


\section{The Pascal adic transformation}

We define the Pascal adic transformation first in terms of its graph,
then we give the cutting and stacking model to which it is isomorphic.

\subsection{The graph construction} \label{sec:pascal-graph}
\begin{figure}[hbt]
\begin{picture}(0,0)%
\includegraphics{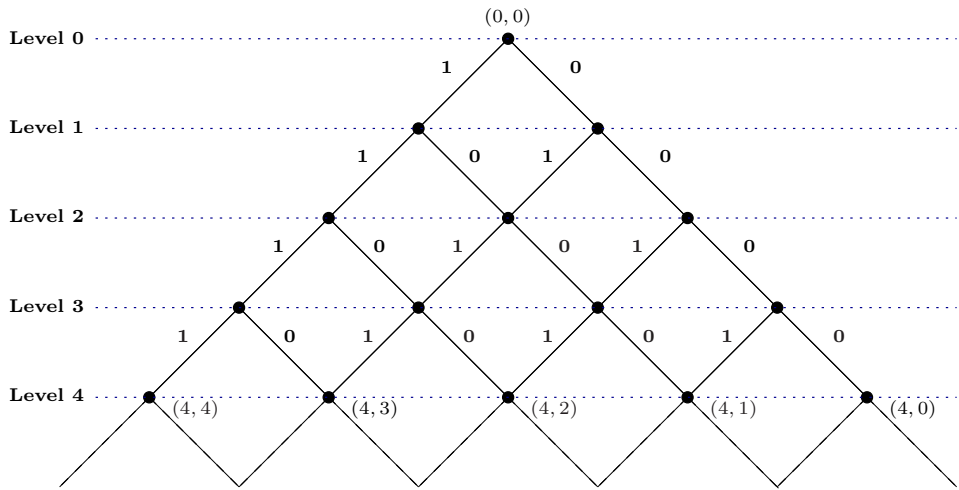}%
\end{picture}%
\setlength{\unitlength}{3947sp}%
\begingroup\makeatletter\ifx\SetFigFont\undefined%
\gdef\SetFigFont#1#2#3#4#5{%
  \reset@font\fontsize{#1}{#2pt}%
  \fontfamily{#3}\fontseries{#4}\fontshape{#5}%
  \selectfont}%
\fi\endgroup%
\begin{picture}(5962,3038)(2451,-3634)
\put(5971,-1009){\makebox(0,0)[lb]{\smash{\SetFigFont{7}{8.4}{\rmdefault}{\bfdefault}{\updefault}{\color[rgb]{0,0,0}0}%
}}}
\put(5161,-1009){\makebox(0,0)[lb]{\smash{\SetFigFont{7}{8.4}{\rmdefault}{\bfdefault}{\updefault}{\color[rgb]{0,0,0}1}%
}}}
\put(4633,-1572){\makebox(0,0)[lb]{\smash{\SetFigFont{7}{8.4}{\rmdefault}{\bfdefault}{\updefault}{\color[rgb]{0,0,0}1}%
}}}
\put(4105,-2135){\makebox(0,0)[lb]{\smash{\SetFigFont{7}{8.4}{\rmdefault}{\bfdefault}{\updefault}{\color[rgb]{0,0,0}1}%
}}}
\put(3506,-2699){\makebox(0,0)[lb]{\smash{\SetFigFont{7}{8.4}{\rmdefault}{\bfdefault}{\updefault}1}}}
\put(5795,-1572){\makebox(0,0)[lb]{\smash{\SetFigFont{7}{8.4}{\rmdefault}{\bfdefault}{\updefault}{\color[rgb]{0,0,0}1}%
}}}
\put(6359,-2135){\makebox(0,0)[lb]{\smash{\SetFigFont{7}{8.4}{\rmdefault}{\bfdefault}{\updefault}1}}}
\put(6922,-2699){\makebox(0,0)[lb]{\smash{\SetFigFont{7}{8.4}{\rmdefault}{\bfdefault}{\updefault}1}}}
\put(5232,-2135){\makebox(0,0)[lb]{\smash{\SetFigFont{7}{8.4}{\rmdefault}{\bfdefault}{\updefault}1}}}
\put(4668,-2699){\makebox(0,0)[lb]{\smash{\SetFigFont{7}{8.4}{\rmdefault}{\bfdefault}{\updefault}1}}}
\put(5795,-2699){\makebox(0,0)[lb]{\smash{\SetFigFont{7}{8.4}{\rmdefault}{\bfdefault}{\updefault}1}}}
\put(5337,-1572){\makebox(0,0)[lb]{\smash{\SetFigFont{7}{8.4}{\rmdefault}{\bfdefault}{\updefault}{\color[rgb]{0,0,0}0}%
}}}
\put(4739,-2135){\makebox(0,0)[lb]{\smash{\SetFigFont{7}{8.4}{\rmdefault}{\bfdefault}{\updefault}{\color[rgb]{0,0,0}0}%
}}}
\put(4175,-2699){\makebox(0,0)[lb]{\smash{\SetFigFont{7}{8.4}{\rmdefault}{\bfdefault}{\updefault}{\color[rgb]{0,0,0}0}%
}}}
\put(5302,-2699){\makebox(0,0)[lb]{\smash{\SetFigFont{7}{8.4}{\rmdefault}{\bfdefault}{\updefault}0}}}
\put(6429,-2699){\makebox(0,0)[lb]{\smash{\SetFigFont{7}{8.4}{\rmdefault}{\bfdefault}{\updefault}0}}}
\put(5901,-2135){\makebox(0,0)[lb]{\smash{\SetFigFont{7}{8.4}{\rmdefault}{\bfdefault}{\updefault}0}}}
\put(6535,-1572){\makebox(0,0)[lb]{\smash{\SetFigFont{7}{8.4}{\rmdefault}{\bfdefault}{\updefault}{\color[rgb]{0,0,0}0}%
}}}
\put(7063,-2135){\makebox(0,0)[lb]{\smash{\SetFigFont{7}{8.4}{\rmdefault}{\bfdefault}{\updefault}{\color[rgb]{0,0,0}0}%
}}}
\put(7626,-2699){\makebox(0,0)[lb]{\smash{\SetFigFont{7}{8.4}{\rmdefault}{\bfdefault}{\updefault}0}}}
\put(3471,-3156){\makebox(0,0)[lb]{\smash{\SetFigFont{7}{8.4}{\rmdefault}{\bfdefault}{\updefault}$(4,4)$}}}
\put(4598,-3156){\makebox(0,0)[lb]{\smash{\SetFigFont{7}{8.4}{\rmdefault}{\bfdefault}{\updefault}{\color[rgb]{0,0,0}$(4,3)$}%
}}}
\put(5725,-3156){\makebox(0,0)[lb]{\smash{\SetFigFont{7}{8.4}{\rmdefault}{\bfdefault}{\updefault}$(4,2)$}}}
\put(6852,-3156){\makebox(0,0)[lb]{\smash{\SetFigFont{7}{8.4}{\rmdefault}{\bfdefault}{\updefault}$(4,1)$}}}
\put(7978,-3156){\makebox(0,0)[lb]{\smash{\SetFigFont{7}{8.4}{\rmdefault}{\bfdefault}{\updefault}{\color[rgb]{0,0,0}$(4,0)$}%
}}}
\put(5433,-692){\makebox(0,0)[lb]{\smash{\SetFigFont{7}{8.4}{\rmdefault}{\bfdefault}{\updefault}{\color[rgb]{0,0,0}$(0,0)$}%
}}}
\put(2451,-1948){\makebox(0,0)[lb]{\smash{\SetFigFont{7}{8.4}{\rmdefault}{\bfdefault}{\updefault}{\color[rgb]{0,0,0}Level 2}%
}}}
\put(2451,-2512){\makebox(0,0)[lb]{\smash{\SetFigFont{7}{8.4}{\rmdefault}{\bfdefault}{\updefault}{\color[rgb]{0,0,0}Level 3}%
}}}
\put(2451,-3075){\makebox(0,0)[lb]{\smash{\SetFigFont{7}{8.4}{\rmdefault}{\bfdefault}{\updefault}{\color[rgb]{0,0,0}Level 4}%
}}}
\put(2451,-1387){\makebox(0,0)[lb]{\smash{\SetFigFont{7}{8.4}{\rmdefault}{\bfdefault}{\updefault}{\color[rgb]{0,0,0}Level 1}%
}}}
\put(2451,-828){\makebox(0,0)[lb]{\smash{\SetFigFont{7}{8.4}{\rmdefault}{\bfdefault}{\updefault}{\color[rgb]{0,0,0}Level 0}%
}}}
\end{picture}
\caption{The Pascal Graph. The number of finite paths from the root to a vertex $(n,k)$ is given by the binomial coefficient $C(n,k)=n!/[k!(n-k)!]$.} \label{fig:pascal-graph}
\end{figure}

\begin{figure}[hbt]
\begin{picture}(0,0)%
\includegraphics{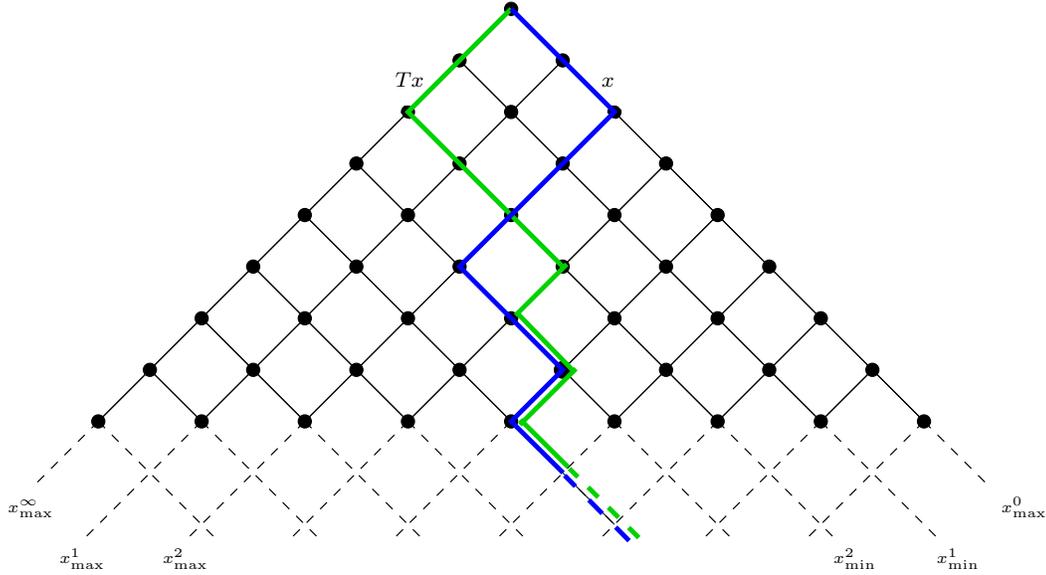}%
\end{picture}%
\setlength{\unitlength}{3947sp}%
\begingroup\makeatletter\ifx\SetFigFont\undefined%
\gdef\SetFigFont#1#2#3#4#5{%
  \reset@font\fontsize{#1}{#2pt}%
  \fontfamily{#3}\fontseries{#4}\fontshape{#5}%
  \selectfont}%
\fi\endgroup%
\begin{picture}(6240,3564)(751,-3207)
\put(6991,-2852){\makebox(0,0)[lb]{\smash{\SetFigFont{6}{7.2}{\familydefault}{\mddefault}{\updefault}{\color[rgb]{0,0,0}$x^0_\text{max}$}%
}}}
\put(6585,-3176){\makebox(0,0)[lb]{\smash{\SetFigFont{6}{7.2}{\familydefault}{\mddefault}{\updefault}{\color[rgb]{0,0,0}$x^1_\text{min}$}%
}}}
\put(5937,-3176){\makebox(0,0)[lb]{\smash{\SetFigFont{6}{7.2}{\familydefault}{\mddefault}{\updefault}{\color[rgb]{0,0,0}$x^2_\text{min}$}%
}}}
\put(751,-2852){\makebox(0,0)[lb]{\smash{\SetFigFont{6}{7.2}{\familydefault}{\mddefault}{\updefault}{\color[rgb]{0,0,0}$x^\infty_\text{max}$}%
}}}
\put(1075,-3176){\makebox(0,0)[lb]{\smash{\SetFigFont{6}{7.2}{\familydefault}{\mddefault}{\updefault}{\color[rgb]{0,0,0}$x^1_\text{max}$}%
}}}
\put(1723,-3176){\makebox(0,0)[lb]{\smash{\SetFigFont{6}{7.2}{\familydefault}{\mddefault}{\updefault}{\color[rgb]{0,0,0}$x^2_\text{max}$}%
}}}
\put(4479,-178){\makebox(0,0)[lb]{\smash{\SetFigFont{8}{9.6}{\rmdefault}{\bfdefault}{\updefault}{\color[rgb]{0,0,0}$x$}%
}}}
\put(3182,-178){\makebox(0,0)[lb]{\smash{\SetFigFont{8}{9.6}{\rmdefault}{\bfdefault}{\updefault}{\color[rgb]{0,0,0}$Tx$}%
}}}
\end{picture}

\caption{$T$ permutes finitely many coordinates and leaves the others unchanged; it ``carries'' a path transversely to the shift transformation. Illustrated in the figure are : $x=0011{\bf 10}0100...$ and $Tx=1100{\bf 01}0100...$.} \label{fig:pascal_image}
\end{figure}

The \emph{Pascal graph} is an infinite planar graph divided into levels $n=0,1,\dots$, with a root vertex at level 0 labeled $(0,0)$, and $n+1$ vertices at each level $n$ labeled $(n,k)$ for $k=0,\dots,n$. From each vertex $(n,k)$ leave two edges; one goes to $(n+1,k+1)$ and is labeled by 1, and the other goes to $(n+1,k)$ and is labeled by 0 --- see Figure \ref{fig:pascal-graph}.
The space $X$ considered is the set of infinite paths going from the
root down the graph, i.e. the set of all $(n,k_n)_{n\ge 1}$, where
$0\le k_n \le n$ and $k_{n+1}=k_n$ or $k_n+1$. The labeling of the
edges produces a natural one-to-one correspondence between $X$ and the
set $\{0,1\}^{\mathbb N}$ of infinite sequences of 0's and 1's. The space $X$
is compact for the product topology, and we denote by $\B$ the Borel
$\sigma$-algebra. Let $d$ be the usual metric on the space
$\{0,1\}^\N$ (letting $d(x,y)=(n+1)^{-1}$ whenever $x$ and $y$
disagree for the first time below the $n$'th level). A cylinder set in
$X$ is a set of the type $\{x\in X\,:\,
x_{i_1}=a_1,x_{i_2}=a_2,\dots,x_{i_s}=a_s\}$, and the family of
cylinder sets generate $\B$. For convenience we will often denote by
$[a_1a_2\dots a_s]$ the cylinder set $\{x\in X\,:\,
x_{1}=a_1,x_{2}=a_2,\dots,x_{s}=a_s\}$. We will refer to a point $x\in
X$ by $(n,k_n(x))_{n\ge 1}$ or by $x_1x_2x_3\dots$, where
$x_1,x_2,\dots$ are the successive labels of the edges of $x$ and
$k_n(x)=\sum_{i=1}^{n}x_i$. We put a partial order on $X$, writing
$x<y$ for $x,y\in X$, whenever $x$ and $y$ coincide below a certain
level $n$ and $x_n<y_n$. In other words, $x$ and $y$ are comparable
with respect to this partial order if for some $n\in\N$ $n_j=y_j$ for
all $j>n$ and there is a permutation $\pi\in\mathcal{S}_n$ such that
$\pi(x_1,\dots,x_n)=(y_1,\dots,y_n)$. Let $X_{\text{min}}$ and
$X_{\text{max}}$ be respectively the set of minimal and maximal
paths. We have
\begin{align*}X_{\text{max}}&=\{x^i_{\text{max}}=0^i1^{\infty}\,:\,i\ge
1\}\cup\{x^0_\text{max}=0^\infty,x^\infty_\text{max}=1^\infty\}\\
X_{\text{min}}&=\{x^i_{\text{min}}=1^i0^{\infty}\,:\,i\ge
1\}.\end{align*}
\begin{definition}
The \emph{Pascal adic transformation} is defined from $X\setminus
X_\text{max}$ onto $X\setminus X_\text{min}$ by $Tx$ = smallest y
greater than $x$.  For every $x\in X\setminus X_\text{max}$ there
are positive integers $n,m$, and $x'\in\{0,1\}^\N$ such that
$x=0^n1^m10x'$. Hence an equivalent definition (illustrated in
Figure \ref{fig:pascal_image}) of $T$ is
\[\boldsymbol{T(0^n1^m10x')=1^m0^n01x'}.\]
\end{definition}

Note that if $x=1^k0^{n-k}\dots$, i.e. $x$ coincides with the minimal
path through the vertex $(n,k)$, then $T^ix$ for $i=1,\dots,C(n,k)-1$
goes through all $C(n,k)$ finite paths from $(0,0)$ to $(n,k)$.

There is a natural way to
extend $T$ bijectively on the whole space $X$ by sending maximal paths
to minimal ones:
\begin{align*}
Tx^i_{\text{max}}&:=x^{i}_\text{min} \quad \text{for all } i\ge 1, \\
Tx^0_\text{max}&:=x^0_\text{max}, \\
Tx^\infty_\text{max}&:=x^\infty_\text{max}.
\end{align*}
Unfortunately, this extension is not continuous at the points
$x^i_\text{max}$,
$x^\infty_\text{max}$.

A simple observation is that if a path $x$ has a ``kink'' at level $n$, i.e. if $x_{n+1}x_{n+2}=10$ --- see Figure \ref{fig:kink} --- then $x$ comes back close to itself after $C(n,k_n(x))$ steps:
\begin{lemma}[The ``Kink'' Lemma]\label{lem:kink_lemma}
If $x\in X$ be a path such that $x_{n+1}x_{n+2}=10$, then $T^{C(n,k_n(x))}x$ and $x$ coincide along the first $n$ coordinates.  
\end{lemma}
\begin{proof}
Consider the following two paths (depicted in Figure \ref{fig:kink})
lying in the orbit of $x$:
\begin{align*}x^+&=0^{n-k_n(x)}1^{k_n(x)}10x_{n+3}x_{n+4}\dots\\x^-&=1^{k_n(x)}0^{n-k_n(x)}01x_{n+3}x_{n+4}\dots.\end{align*}
The adic takes the path $x$ to $x^+$ after $l$ iterations, for a certain integer $l$.  Then the adic applied once to $x^+$ brings us to $x^-$. And finally, it takes $m$ iterations of the adic, for a certain integer $m$, to match up the first $n$ coordinates of $x$.  The total number of iterations we have made is simply the total number of finite paths from $(0,0)$ to $(n,k_n(x))$, so that $l+1+m=C(n,k_n(x))$, establishing the lemma. 

\begin{figure}[hbt]
\begin{picture}(0,0)%
\includegraphics{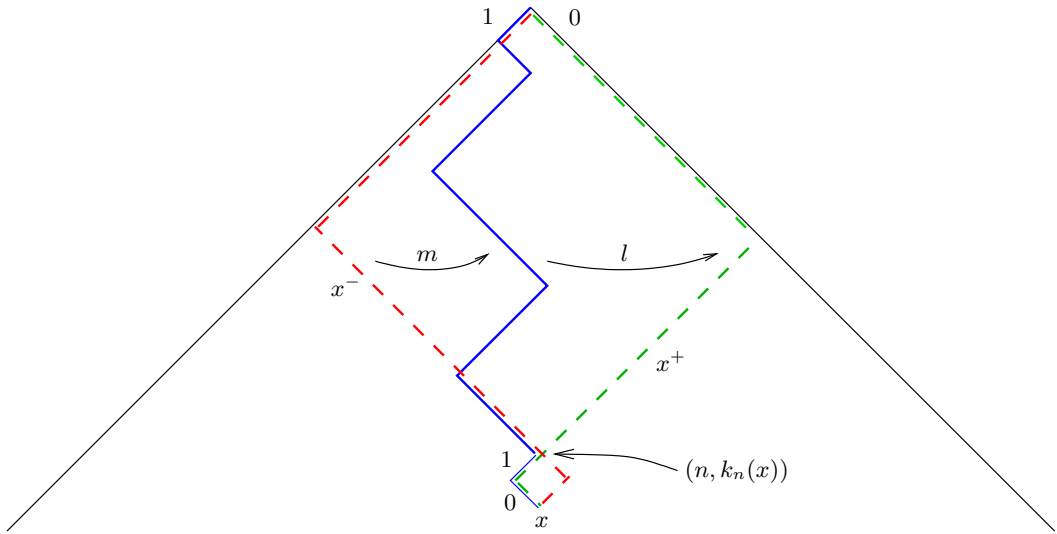}%
\end{picture}%
\setlength{\unitlength}{3947sp}%
\begingroup\makeatletter\ifx\SetFigFont\undefined%
\gdef\SetFigFont#1#2#3#4#5{%
  \reset@font\fontsize{#1}{#2pt}%
  \fontfamily{#3}\fontseries{#4}\fontshape{#5}%
  \selectfont}%
\fi\endgroup%
\begin{picture}(6604,3327)(2573,-4270)
\put(5901,-4230){\makebox(0,0)[lb]{\smash{\SetFigFont{9}{10.8}{\rmdefault}{\bfdefault}{\updefault}{\color[rgb]{0,0,0}$x$}%
}}}
\put(6663,-3248){\makebox(0,0)[lb]{\smash{\SetFigFont{9}{10.8}{\rmdefault}{\bfdefault}{\updefault}{\color[rgb]{0,0,0}$x^+$}%
}}}
\put(4618,-2776){\makebox(0,0)[lb]{\smash{\SetFigFont{9}{10.8}{\rmdefault}{\bfdefault}{\updefault}{\color[rgb]{0,0,0}$x^-$}%
}}}
\put(5567,-1068){\makebox(0,0)[lb]{\smash{\SetFigFont{9}{10.8}{\rmdefault}{\bfdefault}{\updefault}{\color[rgb]{0,0,0}$1$}%
}}}
\put(5156,-2560){\makebox(0,0)[lb]{\smash{\SetFigFont{9}{10.8}{\rmdefault}{\bfdefault}{\updefault}{\color[rgb]{0,0,0}$m$}%
}}}
\put(6441,-2560){\makebox(0,0)[lb]{\smash{\SetFigFont{9}{10.8}{\rmdefault}{\bfdefault}{\updefault}{\color[rgb]{0,0,0}$l$}%
}}}
\put(6114,-1074){\makebox(0,0)[lb]{\smash{\SetFigFont{9}{10.8}{\rmdefault}{\bfdefault}{\updefault}{\color[rgb]{0,0,0}$0$}%
}}}
\put(5707,-4122){\makebox(0,0)[lb]{\smash{\SetFigFont{9}{10.8}{\rmdefault}{\bfdefault}{\updefault}{\color[rgb]{0,0,0}$0$}%
}}}
\put(5685,-3853){\makebox(0,0)[lb]{\smash{\SetFigFont{9}{10.8}{\rmdefault}{\bfdefault}{\updefault}{\color[rgb]{0,0,0}$1$}%
}}}
\put(6843,-3913){\makebox(0,0)[lb]{\smash{\SetFigFont{9}{10.8}{\rmdefault}{\bfdefault}{\updefault}{\color[rgb]{0,0,0}$(n,k_n(x))$}%
}}}
\end{picture}
\caption{The ``Kink'' Lemma} \label{fig:kink}
\end{figure}
\end{proof}

Denote the orbit of a point $x\in X$ by
$\mathcal{O}(x)=\{T^nx\,:\,n\in\Z\}$. $(X,T)$ is not quite a minimal
topological dynamical system (in the sense of a homeomorphism between
compact spaces), but
the Kink Lemma implies that if $10$ appears infinitely many times in
$x$, then $x$ has a dense orbit:
\begin{proposition} $T:X\setminus X_\text{max} \to X\setminus X_\text{min}$ is a homeomorphism, and for every $x\in X$, exactly one of the following holds:
\ben \item  $\mathcal{O}(x)=\{x^0_{\text{max}}\}$ or $\mathcal{O}(x)=\{x^\infty_{\text{max}}\}$ ($x$ is a fixed point) \item there exists $n\ge 1$ such that $x^n_\text{max}\in\mathcal{O}(x)$ (the orbit of $x$ is infinite but not dense)  \item $\overline{\mathcal{O}(x)}=X$ ($x$ has a dense orbit). \een
\end{proposition}

\subsection{Ergodic measures}

If $\tilde{p}=(p,1-p)$ is the probability on $\{0,1\}$ which gives mass
$p$ to $0$ and mass $1-p$ to $1$, then the product measure
$\mu_p=\tilde{p}^{\otimes\N}$ on $\{0,1\}^\N$ is called a Bernoulli
measure and is often denoted by $\B(p,1-p)$.

As noted in \cite{santiago,gibbs} the result that the invariant
ergodic Borel probability measures for the Pascal adic are the
Bernoulli measures $\B(p,1-p)$ is well known; it has been proved by
using the Ergodic Theorem or the Martingale Convergence Theorem. A
more geometric and calculation-free approach developed in \cite{xman2}
permits extension of these results to a wider class of systems, the
generalized Pascal adics. The statement that any
$T$-invariant ergodic measure is a Bernoulli measure can be attributed
to de Finetti in the context of exchangeable processes. The converse,
stating that every Bernoulli measure is ergodic for the Pascal adic,
follows from a result by Hajian, Ito and Kakutani on a system
isomorphic to the Pascal adic defined by interval splitting
\cite{hik}. The connection with adics was made by Vershik
\cite{vershik9}.

In fact the Pascal adic is totally ergodic (every power $T^n$ is
ergodic). This is equivalent to saying that $T$ does not have any
eigenvalues (other than 1) which are roots of unity, which follows
from
the self-similar structure of Pascal's triangle modulo any prime
(a consequence of a result of Lucas \cite{lucas,bollinger}).

\subsection{The cutting and stacking equivalent}

\begin{figure}[hbt]
\includegraphics[width=7cm]{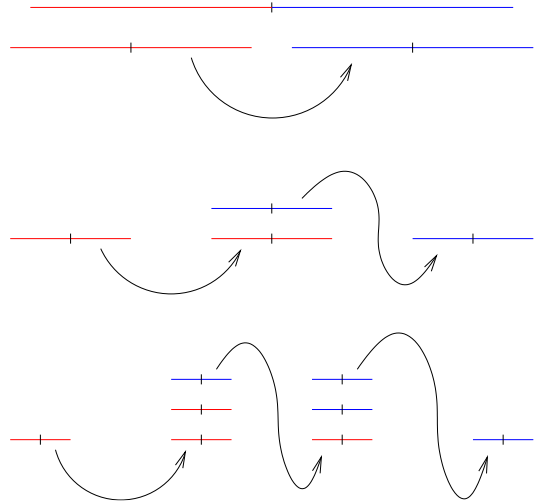}
\caption{The binomial transformation, or cutting and stacking construction of the Pascal adic} \label{pascal_cut_and_stack}
\end{figure}

Start by dividing the unit interval into two equal pieces. At each
step, the stacks are divided into two equal halves, and the right half
of each stack is placed on the bottom of the left half of the
following stack --- see Figure \ref{pascal_cut_and_stack}. If we
repeat indefinitely, the resulting map $T_\text{b}$ (which maps every
open interval of each stack to the one above it) is defined everywhere
except at the dyadic rationals (which correspond to the paths which
are eventually diagonal in the graph construction). Denote by ${m}$
Lebesgue measure, and let $\B([0,1])$ be the $\sigma$-algebra of Borel
sets in $[0,1]$. $([0,1],\B([0,1]),T_\text{b},{m})$ is a
measure-preserving system which we will refer to as the \emph{binomial
transformation}.\\

\begin{proposition}
The systems $([0,1],\B([0,1]),T_\text{b},{m})$ and $(X,\B,T,\mu_{1/2})$ are isomorphic.
\end{proposition}

\begin{proof}
The isomorphism $\psi:[0,1]\setminus\{\text{dyadic rationals}\}\to X\setminus\{y\,:\, x\in X_{\text{min}}\cup X_{\text{max}} \text{ and } y\in \mathcal{O}(x)\}$ is defined by $\psi(\sum_{i=1}^\infty x_i2^{-i})=x_1x_2\dots$. By induction one can check that at step $n$ the bottom and top levels of the $i$'th stack in the cutting and stacking are mapped  respectively to the cylinders $[1^{n-i+1}0^{i-1}]$ and $[0^{i-1}1^{n-i+1}]$. Note also that two points exactly above one another have dyadic expansions which eventually coincide. If $x$ is not a dyadic rational then for some level $n$ it belongs to the right half of the top level in some stack, i.e. $\psi(x)=0^{i-1}1^{n-i+1}0x_{n+1}x_{n+2}\dots$ for some $i\in\{0,\dots,n-1\}$. Therefore $\psi(T_\text{b}(x))=1^{n-i}0^{i}1x_{n+1}x_{n+2}\dots=T(\psi(x))$. Furthermore, $\psi$ takes the Lebesgue measure ${m}$ to the Bernoulli measure $\mu_{1/2}$, since the measures of the levels at step $n$ are equal to $2^{-n}$.
\end{proof}

\begin{remark}
The inverse of this isomorphism also carries any Bernoulli measure $\B(p,1-p)$ on $\{0,1\}^\N$ to the Cantor measure ${m}_p$ on $[0,1]$. Another point of view is to cut and stack with proportions $p$ and $1-p$; then Lebesgue measure carries to the Bernoulli measure $\B(p,1-p)$. \end{remark}

\begin{figure}[hbt]
\includegraphics[width=8cm]{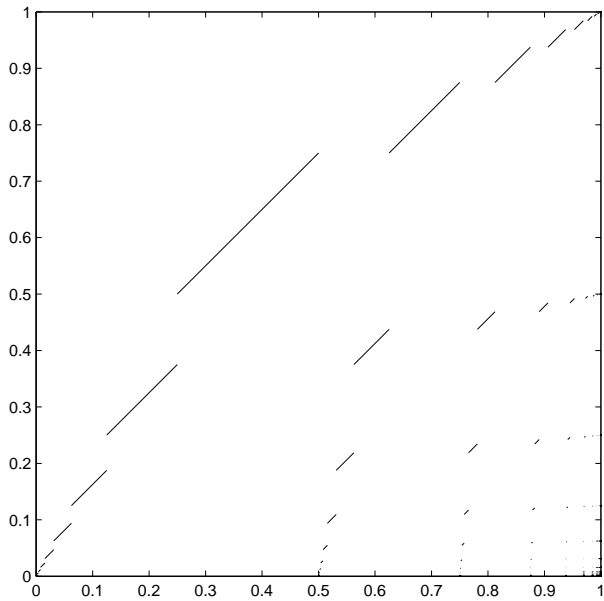}
\caption{The plot of the Pascal adic (or Binomial transformation).} \label{fig:theplot}
\end{figure}

We used \emph{Matlab} to produce the plot of the binomial transformation --- see Figure \ref{fig:theplot}. Note the symmetry with respect to the line $y=1-x$, which can be stated as \[T_\text{b}(1-T_\text{b}x)=1-x.\] In the adic point of view this is equivalent to \[TSTx=Sx,\] where $S$ is the transformation of $\{0,1\}^\N$ which interchanges 0's and 1's. We state this observation as follows:

\begin{proposition}
The Pascal adic $T$ is conjugate to its inverse via the map $S$ which interchanges symbols: $TS=ST^{-1}$.

\end{proposition}

\begin{proof}
Let $x\in X\setminus X_\text{max}$. We can write $x$ in the form $x=0^i1^j\mathbf{10}x'$, where $i,j\ge 0$ and $x'$ is an infinite string of 0's and 1's. Then \begin{align*}
TST(x)&=TS(1^j0^i\mathbf{01}x')\\
&=T(0^j1^i\mathbf{10}S(x'))\\
&=1^i0^j\mathbf{01}S(x')=S(x).\end{align*}
\end{proof}


\section{Limit laws for return times into cylinders}

Considerable attention has been devoted recently to determining the
asymptotic laws of return times or hitting times to ``shrinking
targets'': see, for example,
\cite{pitskel,collet,coelho,paccaut,saussol1,ACS,lacroix,durand-maass} and the references that
they contain.  We establish for the Pascal adic that the limit laws of
return times into typical cylinders, when properly scaled, are
piecewise constant.

For simplicity we assume that $\mu$ is the Bernoulli measure
$\B(1/2,1/2)$, but the steps below can be adapted to the general case.
Fix a generic point $\omega\in X$. Denote by $U_n$ the cylinder
generated by the coordinates $\omega_1, \omega_2,\dots,\omega_{n}$ (as
above we use the notation $U_n=[\omega_1 \omega_2\dots \omega_{n}]$).
Let $\tau_n$ be the first return time (or entrance time) to $U_n$,
i.e. $$\tau_n(x)=\inf\{k\ge 1\,:\, T^k x\in U_n\}.$$ We are interested
in the asymptotics of the return times $\tau_n$, the question being
what is the limit of the following function when the right scaling
$c_n$ is chosen: \[G_n(t)=\fr{1}{\mu(U_n)}\mu\{ x\in U_n\,:\,
c_n\mu(U_n)\tau_{U_n}(x)>t\}.\]

\begin{figure}[hbt]
\begin{picture}(0,0)%
\includegraphics{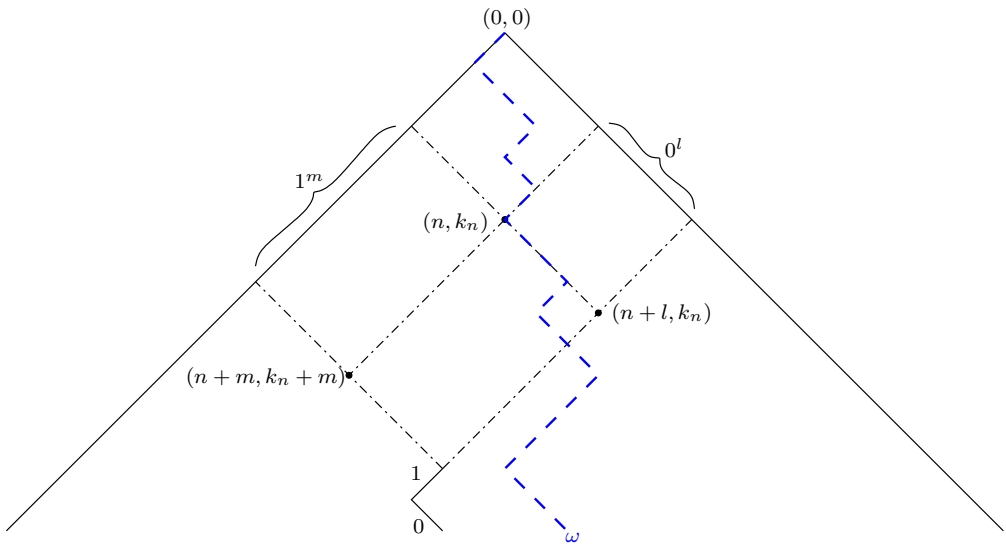}%
\end{picture}%
\setlength{\unitlength}{3947sp}%
\begingroup\makeatletter\ifx\SetFigFont\undefined%
\gdef\SetFigFont#1#2#3#4#5{%
  \reset@font\fontsize{#1}{#2pt}%
  \fontfamily{#3}\fontseries{#4}\fontshape{#5}%
  \selectfont}%
\fi\endgroup%
\begin{picture}(6286,3391)(1184,-2914)
\put(3006,-656){\makebox(0,0)[lb]{\smash{\SetFigFont{8}{9.6}{\rmdefault}{\bfdefault}{\updefault}{\color[rgb]{0,0,0}$1^m$}%
}}}
\put(3750,-2819){\makebox(0,0)[lb]{\smash{\SetFigFont{8}{9.6}{\familydefault}{\mddefault}{\updefault}{\color[rgb]{0,0,0}0}%
}}}
\put(5332,-457){\makebox(0,0)[lb]{\smash{\SetFigFont{8}{9.6}{\rmdefault}{\bfdefault}{\updefault}{\color[rgb]{0,0,0}$0^l$}%
}}}
\put(3735,-2491){\makebox(0,0)[lb]{\smash{\SetFigFont{8}{9.6}{\familydefault}{\mddefault}{\updefault}{\color[rgb]{0,0,0}1}%
}}}
\put(4185,381){\makebox(0,0)[lb]{\smash{\SetFigFont{8}{9.6}{\rmdefault}{\bfdefault}{\updefault}{\color[rgb]{0,0,0}$(0,0)$}%
}}}
\put(4996,-1473){\makebox(0,0)[lb]{\smash{\SetFigFont{8}{9.6}{\rmdefault}{\bfdefault}{\updefault}{\color[rgb]{0,0,0}$(n+l,k_n)$}%
}}}
\put(4708,-2878){\makebox(0,0)[lb]{\smash{\SetFigFont{8}{9.6}{\rmdefault}{\bfdefault}{\updefault}{\color[rgb]{0,0,1}$\omega$}%
}}}
\put(2322,-1878){\makebox(0,0)[lb]{\smash{\SetFigFont{8}{9.6}{\rmdefault}{\bfdefault}{\updefault}{\color[rgb]{0,0,0}$(n+m,k_n+m)$}%
}}}
\put(3810,-901){\makebox(0,0)[lb]{\smash{\SetFigFont{8}{9.6}{\rmdefault}{\bfdefault}{\updefault}{\color[rgb]{0,0,0}$(n,k_n)$}%
}}}
\end{picture}
\caption{The first time an element of $C_{l,m}^n$ returns to $U_n$ is when it enters the cylinder $[\omega_1\dots\omega_n 1^m0^l01]$.}\label{return-time}
\end{figure}

Consider the cylinders $C_{l,m}^n=[\omega_1\dots \omega_{n} 0^l 1^m 1 0]$. Then $U_n=\bigcup_{l,m\ge0}C_{l,m}^n$ (disjoint union up to a set of $\mu$-measure zero).  For $x\in C_{l,m}^n$ and $k_n=\sum_{j=1}^{n}\omega_j$ ---see Figure \ref{return-time}--- the first return time to $U_n$ is given by the following formula:
\begin{proposition} $\tau_{n}(x)=C(n+l,k_n)+C(n+m,k_n+m)-C(n,k_n)$\end{proposition}
\begin{proof}
 If $x\in C_{l,m}^n$, it will return to $U_n$ for the first time when it enters the cylinder $[\omega_1\dots\omega_n 1^m0^l01]$.   After a certain number $N_1$ of iterations $x$ gets mapped to $[0^{n-k_n}1^{k_n}0^l 1^m 1 0]$.  It takes $C(n+l,k_n)-C(n,k_n)$ more iterations to bring $x$ to the cylinder $[0^{n-k_n+l}1^{k_n+m}10]$.  Applying the adic one more time takes us to $[1^{k_n+m}0^{n-k_n+l}01]$.  Similarly, after $C(n+m,k_n+m)-C(n,k_n)$ more iterations we are in $[1^{k_n}0^{n-k_n}1^m 0^l 0 1]$ and after a certain number $N_2$ of iterations $x$ finally re-enters $U_n$.  To conclude note that $N_1+N_2+1=C(n,k_n)$.
\end{proof}

Letting $t_{l,m}^n=\fr{1}{2^n}c_n\tau_{n}$ for $n\geq 1, l,m \geq 0$, we have
\begin{align*}
G_n(t) & =2^n\sum_{l,m\ge 0}\mu\{x\in C_{l,m}^n\,:\, t_{l,m}^n > t\}
       =2^n\sum_{l,m\ge 0\,:\, t_{l,m}^n>t}\mu(C_{l,m}^n) \\
      & =2^n\sum_{l,m\ge 0\,:\, t_{l,m}^n>t}\fr{1}{2^{n+l+m+2}}
       =\fr{1}{4}\sum_{l,m\ge 0\,:\, t_{l,m}^n>t}\fr{1}{2^{l+m}}.
\end{align*}

Using Stirling's Formula, we know that $C(n,k_n)/2^n\approx 1/\sqrt{2\pi(n-k_n)}$, so that when $c_n$ grows faster than $\sqrt{n}$, $\lim_{n\to\infty}t_{l,m}^n=0$ and therefore $\lim_{n\to\infty}G_n(t)=\II_{(-\infty,0]}(t)$; and when $c_n$ grows slower than $\sqrt{n}$, then $\lim_{n\to\infty}t_{l,m}^n=\infty$, which implies that $\lim_{n\to\infty}G_n(t)=\II_{(-\infty,\infty)}$. The interesting scaling is $c_n=\sqrt{n}$. Then using again Stirling's formula and the fact that $k_n/n\to1/2$ $\mu$-a.e. we get
\begin{align*}
t_{0,0}^n &=\fr{1}{2^n}\,\sqrt{n}\,C(n,k_n)=\fr{1}{2^n}\,\sqrt{n}\,\fr{n!}{k_n!(n-k_n)!} \\
          &\approx\fr{1}{2^n}\,\sqrt{n}\,\fr{\fr{n^n}{e^n}\sqrt{2\pi n}}{\fr{k_n^{k_n}}{e^{k_n}}\sqrt{2\pi k_n}\fr{(n-k_n)^{n-k_n}}{e^{n-k_n}}\sqrt{2\pi (n-k_n)}} \\
          &\approx\fr{1}{2^n\sqrt{2\pi}}\,\sqrt{n}\,\fr{n^{n+1/2}}{k_n^{k_n+1/2}(n-k_n)^{n-k_n+1/2}} \\
          &\approx\fr{1}{2^n\sqrt{2\pi}}\,\sqrt{\fr{n}{n-k_n}}\,\left(\fr{n}{k_n}\right)^{k_n+1/2}\left(\fr{n}{n-k_n}\right)^{n-k_n} \\
          &\underset{n\to\infty}{\longrightarrow} \sqrt{\fr{2}{\pi}}.
\end{align*}
\noindent Then
\begin{align*}
t_{l,m}^n &= \fr{1}{2^n}\,\sqrt{n}\,\left[C(n+m,k_n+m)+C(n+l,k_n)-C(n,k_n)\right] \\
          &= t_{0,0}^n\left[\fr{(n+m)\dots(n+1)}{(k_n+m)\dots(k+1)}+\fr{(n+l)\dots(n+1)}{(n+l-k_n)\dots(n+1-k_n)}-1\right] \\
          &\underset{n\to\infty}{\longrightarrow}\sqrt{\fr{2}{\pi}}\;(2^l+2^m-1).
\end{align*}
We can now easily deduce that $G_n$ converges to a step function whose heights are computed below.
Let $t_{i,j}=\sqrt{\fr{2}{\pi}}\;(2^i+2^j-1)$ for $i \geq 0$ and $0
\leq j \leq i$.

\noindent Assume that $t_{i,j}\le t<t_{i,j+1}$ and $j+1\le i$; then

\begin{align*}
\lim_{n\to\infty}G_n(t) &= \fr{1}{4}\lim_{n\to\infty}\sum_{l,m\ge 0 : t_{l,m}^n>t}\fr{1}{2^{l+m}} \\ &=\fr{1}{4}\sum_{l,m\ge 0 : t_{l,m}>t}\fr{1}{2^{l+m}} \\ &=\fr{1}{4}\left(\sum_{l=m\ge 0 : t_{l,m}>t}\fr{1}{2^{l+m}}+2\sum_{l>m\ge 0 : t_{l,m}>t}\fr{1}{2^{l+m}}\right)
\\ &=\fr{1}{4}\left(\sum_{l=i}^{\infty}\fr{1}{2^{2l}}+2\sum_{m=j+1}^{i-1}\fr{1}{2^{i+m}}+2\sum_{l=i+1}^{\infty}\sum_{m=0}^{l-1}\fr{1}{2^{l+m}}\right)
\\ &=2^{-1-2i}(2^{i+1}+2^{i-j}-2).
\end{align*}

Similarly, if $t_{i,i}\le t<t_{i+1,0}$ (the case $i=j$), then
\begin{align*}
\lim_{n\to\infty}G_n(t) &= \fr{1}{4}\left(\sum_{l=i+1}^{\infty}\fr{1}{2^{2l}}+2\sum_{l=i+1}^{\infty}\sum_{m=0}^{l-1}\fr{1}{2^{l+m}}\right)\\
&= 2^{-2(i+1)}(2^{i+2}-1).
\end{align*}

To summarize:
\begin{thm}\label{thm:returntimes}
$G_n(t)=\fr{1}{\mu(U_n)}\,\mu\{ x\in U_n\,:\, \sqrt{n}\mu(U_n)\tau_{n}(x)>t\}$ converges pointwise as $n \to \infty$ to the piecewise constant function
\begin{center}
  $\begin{cases}
    1 & \text{ if }  t<0 \\
    2^{-1-2i}(2^{i+1}+2^{i-j}-2) & \text{ if } t_{i,j}\le t <t_{i,j+1} \text{ and } i>j+1   \\
    2^{-2(i+1)}(2^{i+2}-1) & \text{ if } t_{i,i}\le t<t_{i+1,0},
\end{cases}$
\end{center}
where $t_{i,j}=\sqrt{\fr{2}{\pi}}\,(2^i+2^j-1)$ for $i \geq 0$ and $0
\leq j \leq i$.
\end{thm}
\begin{remark}
If
$p\ne 1/2$ and $\mu$ is the Bernoulli measure $\B(p, 1-p)$, then $G_n(t)$ defined as above still converges to a piecewise constant function, but the formula for the function is rather messy.
\end{remark}


\section{A countable-substitution subshift}

Stationary adics are topologically conjugate to substitution or odometer symbolic dynamical systems \cite{forrest,livshitz,durand-host-skau}. Here we use the idea described in \cite{durand-host-skau,host1} to associate a substitution to an adic to show that the Pascal adic is isomorphic to a subshift whose language is determined by countably many substitutions. For basic definitions and general information about substitution systems see \cite{queffelec,fogg}. We will use the alphabet $\{ a,b\}$ for our subshift to decrease confusion about the various representations of the Pascal adic.

\begin{figure}[ht]
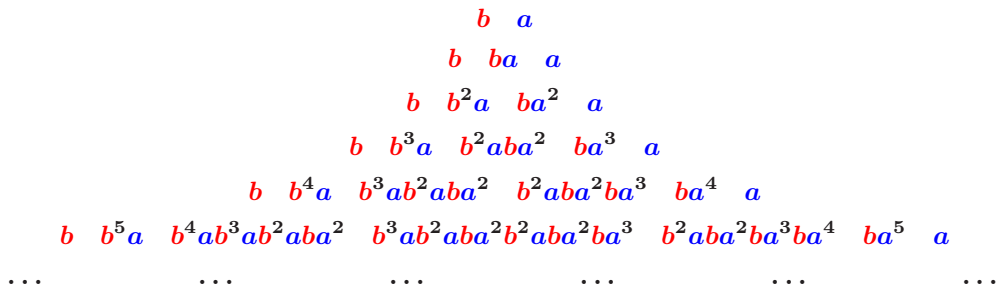

\begin{gather*}
\boldsymbol{\textcolor{red}{b} \quad \textcolor{blue}{a}} \\
\boldsymbol{\textcolor{red}{b} \quad \textcolor{red}{b}\textcolor{blue}{a} \quad \textcolor{blue}{a}} \\
\boldsymbol{\textcolor{red}{b} \quad \textcolor{red}{b}^2\textcolor{blue}{a} \quad \textcolor{red}{b}\textcolor{blue}{a}^2 \quad \textcolor{blue}{a}}\\
\boldsymbol{\textcolor{red}{b} \quad \textcolor{red}{b}^3\textcolor{blue}{a} \quad \textcolor{red}{b}^2\textcolor{blue}{a}\textcolor{red}{b}\textcolor{blue}{a}^2 \quad \textcolor{red}{b}\textcolor{blue}{a}^3 \quad \textcolor{blue}{a}}\\
\boldsymbol{\textcolor{red}{b} \quad \textcolor{red}{b}^4\textcolor{blue}{a} \quad \textcolor{red}{b}^3\textcolor{blue}{a}\textcolor{red}{b}^2\textcolor{blue}{a}\textcolor{red}{b}\textcolor{blue}{a}^2 \quad \textcolor{red}{b}^2\textcolor{blue}{a}\textcolor{red}{b}\textcolor{blue}{a}^2\textcolor{red}{b}\textcolor{blue}{a}^3 \quad \textcolor{red}{b}\textcolor{blue}{a}^4 \quad \textcolor{blue}{a}}\\
\boldsymbol{\textcolor{red}{b} \quad \textcolor{red}{b}^5\textcolor{blue}{a} \quad \textcolor{red}{b}^4\textcolor{blue}{a}\textcolor{red}{b}^3\textcolor{blue}{a}\textcolor{red}{b}^2\textcolor{blue}{a}\textcolor{red}{b}\textcolor{blue}{a}^2 \quad \textcolor{red}{b}^3\textcolor{blue}{a}\textcolor{red}{b}^2\textcolor{blue}{a}\textcolor{red}{b}\textcolor{blue}{a}^2\textcolor{red}{b}^2\textcolor{blue}{a}\textcolor{red}{b}\textcolor{blue}{a}^2\textcolor{red}{b}\textcolor{blue}{a}^3 \quad \textcolor{red}{b}^2\textcolor{blue}{a}\textcolor{red}{b}\textcolor{blue}{a}^2\textcolor{red}{b}\textcolor{blue}{a}^3\textcolor{red}{b}\textcolor{blue}{a}^4 \quad \textcolor{red}{b}\textcolor{blue}{a}^5 \quad  \textcolor{blue}{a}}\\
\boldsymbol{\dots\hspace{2cm}\dots\hspace{2cm}\dots\hspace{2cm}\dots\hspace{2cm}\dots\hspace{2cm}\dots}
\end{gather*}
\caption{``The Pascal triangle of words'', showing the construction of the words by successive concatenations from level 1 to level 6. The word written at the ``vertex'' $(n,k)$ of this triangle represents the coding of the cylinder $[1^k0^{n-k}]$.}
\label{fig:blocks}
\end{figure}

Consider the substitutions $\zeta_i:\{0,1,\dots,i\}\to\{a,b\}^*=$ all finite words on the alphabet $\{ a,b\}$ defined by
\begin{align*}
\zeta_i(i)&=b \quad \zeta_i(0)=a \quad \text{for all } i\ge 1 \\
\zeta_i(j)&=\zeta_{i-1}(j)\zeta_{i-1}(j-1)  \quad \text{for all } i\ge 2 \text{ and } 0<j<i.
\end{align*}
Denote by $\mathcal{L}(\zeta_i)$ the language associated to $\zeta_i$,
by which we mean the subset of $\{ a,b\}^*$ consisting of all subwords
of all the $\zeta_i(j), i \geq 1, 0 \leq j \leq i$. Let $\Sigma$ be
the subshift consisting of all sequences in $\{a,b\}^\Z$ all of whose
subwords belong to
$\mathcal{L}(\zeta)=\cup_i\mathcal{L}(\zeta_i)$. Figure
\ref{fig:blocks} shows how the ``basic'' words in the language
$\mathcal{L}(\zeta)$ can be obtained by successive concatenations. It is not very hard to see that these basic words give the codings of the cylinders $[1^k0^{n-k}]$ in the
Pascal graph under the action of the Pascal adic transformation
according to the first edge (left $\sim$ label 1 $\sim b$, right $\sim$
label 0 $\sim a$)---or, in the cutting and stacking representation,
the coding according to visits to $[0,1/2) \sim b$ or $[1/2,1] \sim
a$ (see Figure \ref{fig:isomorphism2}). (Just note that the cylinder $[1^k0^{n-k}]$ corresponds to the bottom level of the $k$th stack in the cutting and stacking at stage $n$, thus its coding is the word seen going up the stack (Figure \ref{fig:isomorphism2}), which is the word $\zeta_n(k)$ shown at vertex $(n,k)$ in Figure \ref{fig:blocks}).

\begin{thm}\label{thm:pascal-sub-iso}
Let $(X,T)$ denote the Pascal adic. There are a countable subset
$X'\subset X$ and a one-to-one Borel measurable map $\phi:X\setminus
X'\to\Sigma$ such that $\phi\circ T=\sigma\circ\phi$ on $X\setminus
X'$. Consequently, for all $\alpha\in(0,1)$, for each nonatomic
ergodic measure $\mu_\alpha$ on $(X,T)$, if
$\nu_\alpha=\mu_\alpha\circ\phi^{-1}$, then $(X,T,\mu_\alpha)$ and
$(\Sigma,\sigma,\nu_\alpha)$ are measurably isomorphic.
\end{thm}
\begin{proof}
Let $P=\{[0],[1]\}$, where we denote as previously $[0]=\{x\in X\,:\,x_0=0\}$ and $[1]=\{x\in X\,:\,x_0=1\}.$ Let $X'$ be the countable subset of $X$ consisting  of all paths which are eventually diagonal, i.e. $x\in X'$ if and only if there exists $n\in\N$ such that either $x_k=0$ for all $k\ge n$, or $x_k=1$ for all $k\ge n$. Let $\phi$ be the coding of the Pascal adic $(X,T)$ by the partition $P$, where $[0]$ is coded with $a$ and $[1]$ with $b$. More precisely, for every $x\in X\setminus X'$ we define \[\phi(x)=(\dots,\omega_{-2},\omega_{-1},\omega_0,\omega_1,\omega_2,\dots),\] where \[\omega_i=\begin{cases}a \text{ if } (T^ix)_0=0 \\ b \text{ if } (T^ix)_0=1.\end{cases}\]
\begin{figure}[hbt]
\begin{picture}(0,0)%
\includegraphics{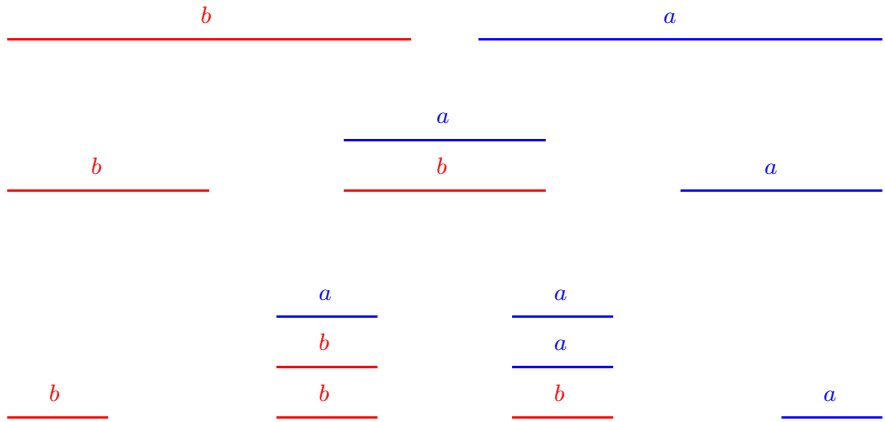}%
\end{picture}%
\setlength{\unitlength}{3947sp}%
\begingroup\makeatletter\ifx\SetFigFont\undefined%
\gdef\SetFigFont#1#2#3#4#5{%
  \reset@font\fontsize{#1}{#2pt}%
  \fontfamily{#3}\fontseries{#4}\fontshape{#5}%
  \selectfont}%
\fi\endgroup%
\begin{picture}(5541,2612)(879,-2671)
\put(2117,-167){\makebox(0,0)[lb]{\smash{\SetFigFont{9}{10.8}{\rmdefault}{\bfdefault}{\updefault}{\color[rgb]{1,0,0}$b$}%
}}}
\put(5024,-167){\makebox(0,0)[lb]{\smash{\SetFigFont{9}{10.8}{\rmdefault}{\bfdefault}{\updefault}{\color[rgb]{0,0,1}$a$}%
}}}
\put(5658,-1118){\makebox(0,0)[lb]{\smash{\SetFigFont{9}{10.8}{\rmdefault}{\bfdefault}{\updefault}{\color[rgb]{0,0,1}$a$}%
}}}
\put(1430,-1118){\makebox(0,0)[lb]{\smash{\SetFigFont{9}{10.8}{\rmdefault}{\bfdefault}{\updefault}{\color[rgb]{1,0,0}$b$}%
}}}
\put(3597,-801){\makebox(0,0)[lb]{\smash{\SetFigFont{9}{10.8}{\rmdefault}{\bfdefault}{\updefault}{\color[rgb]{0,0,1}$a$}%
}}}
\put(3597,-1118){\makebox(0,0)[lb]{\smash{\SetFigFont{9}{10.8}{\rmdefault}{\bfdefault}{\updefault}{\color[rgb]{1,0,0}$b$}%
}}}
\put(1165,-2543){\makebox(0,0)[lb]{\smash{\SetFigFont{9}{10.8}{\rmdefault}{\bfdefault}{\updefault}{\color[rgb]{1,0,0}$b$}%
}}}
\put(2857,-1909){\makebox(0,0)[lb]{\smash{\SetFigFont{9}{10.8}{\rmdefault}{\bfdefault}{\updefault}{\color[rgb]{0,0,1}$a$}%
}}}
\put(4337,-1909){\makebox(0,0)[lb]{\smash{\SetFigFont{9}{10.8}{\rmdefault}{\bfdefault}{\updefault}{\color[rgb]{0,0,1}$a$}%
}}}
\put(4337,-2226){\makebox(0,0)[lb]{\smash{\SetFigFont{9}{10.8}{\rmdefault}{\bfdefault}{\updefault}{\color[rgb]{0,0,1}$a$}%
}}}
\put(6028,-2543){\makebox(0,0)[lb]{\smash{\SetFigFont{9}{10.8}{\rmdefault}{\bfdefault}{\updefault}{\color[rgb]{0,0,1}$a$}%
}}}
\put(2857,-2543){\makebox(0,0)[lb]{\smash{\SetFigFont{9}{10.8}{\rmdefault}{\bfdefault}{\updefault}{\color[rgb]{1,0,0}$b$}%
}}}
\put(2857,-2226){\makebox(0,0)[lb]{\smash{\SetFigFont{9}{10.8}{\rmdefault}{\bfdefault}{\updefault}{\color[rgb]{1,0,0}$b$}%
}}}
\put(4337,-2543){\makebox(0,0)[lb]{\smash{\SetFigFont{9}{10.8}{\rmdefault}{\bfdefault}{\updefault}{\color[rgb]{1,0,0}$b$}%
}}}
\end{picture}
\caption{Coding by $a$'s and $b$'s in the cutting and stacking.}
\label{fig:isomorphism2}
\end{figure}
The map $\phi$ intertwines $T$ and $\sigma$ (this is clear) and establishes a measure-theoretic isomorphism between $(X,T,\mu_{\alpha})$ and $(\Sigma,\sigma,\nu_{\alpha})$. To check this, it is enough to show that $\phi$ is one-to-one, or, equivalently, that $P$ is a generating partition for $T$. Let $x,y\in X$ be two different paths which are not eventually diagonal. Since $x\neq y$ there is a smallest integer $n$ such that $x_n\neq y_n$, and we can assume that $x_n=1$ and $y_n=0$. In addition, since $x\notin X_{\max}$, $x_{n+j}=0$ for some smallest $j\ge 1$. Consequently, the $P$-names of $x$ and $y$ coincide until $x$ and $y$ get mapped after $N$ iterations to the cylinder $[0^l1^m]$ for some $N,l,m\ge 1$ (where $m=k_{n-1}(x)=k_{n-1}(y)$, $l=n-k_{n-1}(x)$, and $0\le N \le C(n,k_{n-1}(x))$). Then, since $T^Nx=0^l1^{m+j+1}0\dots$ and $T^Ny=0^l1^m0\dots$, it follows that the $P$-name of $T^Nx$ is $b^{m+j+1}a\dots$ whereas the $P$-name of $T^Ny$ is $b^ma\dots$, showing that $\phi(x)\neq\phi(y)$ --- see Figure \ref{fig:pascal_sub_isom2}. 
\begin{figure}[hbt]
\begin{picture}(0,0)%
\includegraphics{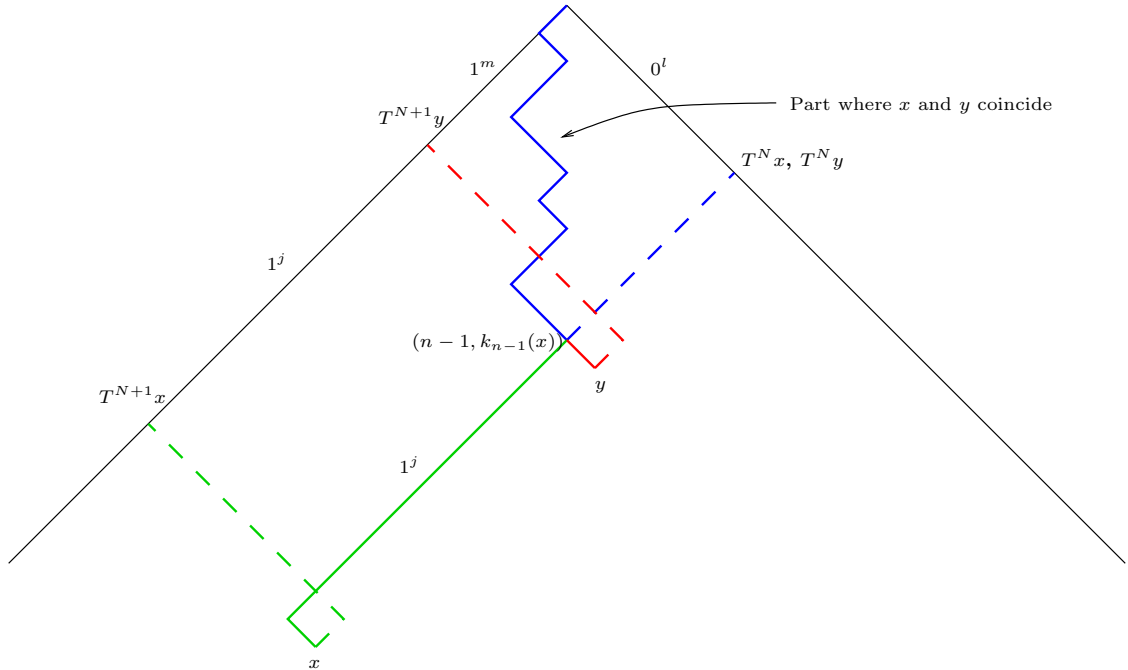}%
\end{picture}%
\setlength{\unitlength}{3947sp}%
\begingroup\makeatletter\ifx\SetFigFont\undefined%
\gdef\SetFigFont#1#2#3#4#5{%
  \reset@font\fontsize{#1}{#2pt}%
  \fontfamily{#3}\fontseries{#4}\fontshape{#5}%
  \selectfont}%
\fi\endgroup%
\begin{picture}(7037,4213)(2365,-4526)
\put(4261,-4499){\makebox(0,0)[lb]{\smash{\SetFigFont{7}{8.4}{\rmdefault}{\bfdefault}{\updefault}{\color[rgb]{0,0,0}$x$}%
}}}
\put(6059,-2746){\makebox(0,0)[lb]{\smash{\SetFigFont{7}{8.4}{\rmdefault}{\bfdefault}{\updefault}{\color[rgb]{0,0,0}$y$}%
}}}
\put(6979,-1343){\makebox(0,0)[lb]{\smash{\SetFigFont{7}{8.4}{\rmdefault}{\bfdefault}{\updefault}{\color[rgb]{0,0,0}$T^Nx$, $T^Ny$}%
}}}
\put(4832,-3272){\makebox(0,0)[lb]{\smash{\SetFigFont{7}{8.4}{\rmdefault}{\bfdefault}{\updefault}{\color[rgb]{0,0,0}$1^j$}%
}}}
\put(5270,-773){\makebox(0,0)[lb]{\smash{\SetFigFont{7}{8.4}{\rmdefault}{\bfdefault}{\updefault}{\color[rgb]{0,0,0}$1^m$}%
}}}
\put(7286,-992){\makebox(0,0)[lb]{\smash{\SetFigFont{7}{8.4}{\familydefault}{\mddefault}{\updefault}{\color[rgb]{0,0,0}Part where $x$ and $y$ coincide}%
}}}
\put(2946,-2834){\makebox(0,0)[lb]{\smash{\SetFigFont{7}{8.4}{\rmdefault}{\bfdefault}{\updefault}{\color[rgb]{0,0,0}$T^{N+1}x$}%
}}}
\put(4700,-1080){\makebox(0,0)[lb]{\smash{\SetFigFont{7}{8.4}{\rmdefault}{\bfdefault}{\updefault}{\color[rgb]{0,0,0}$T^{N+1}y$}%
}}}
\put(3999,-2000){\makebox(0,0)[lb]{\smash{\SetFigFont{7}{8.4}{\rmdefault}{\bfdefault}{\updefault}{\color[rgb]{0,0,0}$1^j$}%
}}}
\put(6410,-773){\makebox(0,0)[lb]{\smash{\SetFigFont{7}{8.4}{\rmdefault}{\bfdefault}{\updefault}{\color[rgb]{0,0,0}$0^l$}%
}}}
\put(4910,-2481){\makebox(0,0)[lb]{\smash{\SetFigFont{7}{8.4}{\familydefault}{\mddefault}{\updefault}{\color[rgb]{0,0,0}$(n-1,k_{n-1}(x))$}%
}}}
\end{picture}
\caption{After $N$ steps $x$ and $y$ are mapped to the ``maximal'' cylinder
  $[0^l1^m]$. At step $N+1$, $x$ is mapped to $1^{m+j-1}0^l$ whereas
  $y$ is mapped to $1^m0^l$. This shows that $x$ will then stay in the
  cylinder $[1]$ coded by $b$ for $m+j-1$ steps until it is mapped to
  the cylinder $[0]$ coded by $a$, whereas it will only take $m$ steps
  for $y$ to get mapped into $[0]$.}
\label{fig:pascal_sub_isom2}
\end{figure}

\end{proof}


\subsection{Complexity}

It is easy to see that the subshift $(\Sigma,\sigma)$ has topological
entropy 0. Finer measures of the size or richness of symbolic
dynamical systems can be drawn from asymptotics of the complexity
function, which for each $n$ gives the number $p_n$ of $n$-blocks
found in sequences in the system. See
\cite{ferenczi3,ferenczi4,troubetzkoy,fogg} for some examples of
results on the complexity functions of various systems. Usually one
finds upper or lower estimates on the growth rate of a complexity
function. Here we show that for the subshift $(\Sigma, \sigma)$ that
we have associated to the Pascal adic, $p_n$ is asymptotic to $n^3/6$.

Let $B_{n,j}=\zeta_j(n)$ denote the block at the vertex $(n,j)$ in
Figure \ref{fig:blocks} (for example $B_{3,2}=b^2a$). We continue to refer to the $B_{n,j}$'s as \emph{basic blocks}. These blocks satisfy the recurrence formula
\begin{equation}\label{recursion}
B_{n,j}=B_{n-1,j}B_{n-1,j-1},\text{ for }0\le j\le n.
\end{equation}
Every basic block at level $n$ can be written as the product (we use interchangeably ``product'' and ``concatenation'') of two basic blocks from level $n-1$, and every basic block at level $n-1$ can be written as the product of two basic blocks from level $n-2$, and so forth... Thus we have a \emph{hierarchical decomposition} of the basic blocks. We can view any block $B_{n,j}$ as the product of basic blocks from level $l$ for any $l\le n$, depending on how far back in the hierarchy we want to look.

\begin{definition}
We say $B\in\mathcal{L}(\zeta)$ is a \emph{new block at level $l$} if $B$ appears as a subblock of one of the basic blocks $B_{l,k}$ at level $l$, but does not appear as a subblock at level $l-1$ (and a fortiori at any lower levels). 
\end{definition}

We would like to count the number of different new $n$-blocks at level $l$ for $l=1,2,\dots$.  Although the recursive construction of the blocks is simple, it is not clear how to count precisely the different $n$-blocks since a given $n$-block can appear many times at the same level.  The concatenation of two consecutive basic blocks at level $l$ will result in the formation of new blocks at level $l+1$, but how can we tell whether a block which overlaps two such basic blocks didn't appear higher in the Pascal triangle of words (see Figure \ref{fig:blocks})? After which level will we have seen all the different $n$-blocks? The following lemma provides an answer to these questions.


\begin{lemma}\label{new-blocks}
Let $B$ be an $n$-subblock at level $l+1$ of $B_{l+1,j}$, for $2\le j\le l-1$. Assume that $l,j$ are such that $|B_{l-1,j-1}|\ge\max\{n-(l-j),n-(j-1)\}$ (which guarantees that $B$ does not overlap $B_{l-1,j}$ or $B_{l-1,j-2}$). Then $B$ is a new $n$-block at level $l+1$ if and only if $B$ contains $a^{l-j}b^{j-1}$.
\end{lemma}

\begin{example}
Observe in Figure \ref{fig:blocks} that the only new 5-blocks which appear at level 6 as subblocks of $B_{6,4},B_{6,3},B_{6,2}$ are \[b{\boldsymbol{ab^3}},{\boldsymbol{ab^3}}a,b{\boldsymbol{a^2b^2}},{\boldsymbol{a^2b^2}}a,b{\boldsymbol{a^3b}},{\boldsymbol{a^3b}}a.\] They all contain the block $\boldsymbol{a^{5-j}b^{j-1}}$ for some $j$.
\end{example}

\begin{proof}
 Suppose that $B$ is a new block at level $l+1$ not containing the subblock $a^{l-j}b^{j-1}$. For example, in the case $B_{l,j}=B_{6,4}$ and $B_{l,j-1}=B_{6,3}$, we have the following picture: \[B_{l,j}= b^4a|\overbrace{b^3ab^2ab\boldsymbol{a^2}}^{B_{l-1,j-1}} \quad \overbrace{\boldsymbol{b^3}ab^2aba^2}^{B_{l-1,j-1}}|b^2aba^2ba^3=B_{l,j-1},\] \[B_{l,j}=b^4a|b^3\underbrace{ab^2ab\boldsymbol{aa \quad bbb}}_{B}\boldsymbol{a}\underbrace{ab^2aba^2|b^2}_{B}aba^2ba^3=B_{l,j-1}.\] \begin{center}(Note: We use ``$|$'' to symbolize where the concatenation at the previous level took place.)\end{center} Observe that $a^{l-j}$ is a right factor of $B_{l,j}$, and that $b^{j-1}$ is a left factor of $B_{l,j-1}$. Consequently, there exists $k$ such that either $B=\omega b^k$ and $1\le k\le j-2$, or $B=a^k\omega$ and $1\le k\le l-j$, where $\omega$ is respectively a right or left factor of $B_{l-1,j-1}$. The latter follows from the hypothesis that $|B_{l-1,j-1}|\ge\max\{n-(l-j),n-(j-1)\}\ge |B|-k$. Suppose that $B=\omega b^k$ (the case  $B=a^k\omega$ is similar). Since $B_{l,j-1}=B_{l-1,j-1}B_{l-1,j-2}=B_{l-1,j-1}b^{j-2}\dots$, it follows that $B$ is a subblock of $B_{l,j-1}$. This contradicts the fact that $B$ is a new block at level $l+1$.

Conversely, observe that for $1\le j\le l-1$, $B_{l,j}=b^ja\dots ba^{l-j}$. Thus the first time the block $a^{l-j}b^{j-1}$ will be seen is at level $l+1$, as a subblock of $B_{l,j}B_{l,j-1}$.
\end{proof}

\noindent As a corollary we get:
\begin{lemma}\label{lem:n+2}
All $n$-blocks are seen at level $n+2$ as subblocks.
\end{lemma}

\begin{proof}
Assume there exists a new $n$-block $B$ at level $n+3$. First note that $B$ cannot be a subblock of any of the following ``edge'' blocks: \[ B_{n+3,n+3}=b,\,B_{n+3,n+2}=b^{n+2}a,\,B_{n+3,1}=ba^{n+2},\,B_{n+3,0}=a.\] Otherwise $B=b^n,b^{n-1}a,a^n,$ or $ba^{n-1}$, and those blocks are already seen at level $n+1$. Therefore $B$ is coming from the concatenation of $B_{n+2,j}$ and $B_{n+2,j-1}$, for some $j$ with $2\le j\le n+1$. By Lemma \ref{new-blocks} (applied in the case $l=n+2$), $B$ must contain the subblock $a^{n+2-j}b^{j-1}$, which is impossible.
\end{proof}

\begin{thm}\label{thm:complexity}
$\ds\lim_{n\to\infty}\frac{p_n}{n^3}=\frac{1}{6}.$
\end{thm}

\begin{proof}
Fix $n$. Consider the Pascal triangle of words from level 1 to level $n+1$. We will count the number of $n$-blocks by estimating how many new $n$-blocks are created from one level to the next by concatenation of two adjacent basic blocks. We divide the triangle into two disjoint regions: let \emph{region I} be the subset of the triangle formed by all blocks up to level $\sqrt{2n}$, as well as all blocks below the level $\sqrt{2n}$ located at the vertices $(l,k)$, for $\sqrt{2n}\le l\le n+1$ and $k=0,1,l-1,l$. Let \emph{region II} be the complement of region I in the triangle considered --- see Figure \ref{fig:block-count}. Let $p_\text{I}$ be the sum for $l=1$ to $n+1$ of the number of new $n$-blocks at level $l$ appearing as subblocks of the product of two adjacent basic blocks, one of them at least belonging to region $\text{I}$. Similarly, let $p_\text{II}$ denote the sum for $\sqrt{2n}\le l\le n+1$ of the number of new $n$-blocks at level $l$ appearing as subblocks of the product of two consecutive basic blocks, both of them belonging to region $\text{II}$. By Lemma \ref{lem:n+2}, since all $n$-blocks are seen at level $n+2$, we have the following inequality: \[p_{\text{II}} \le p_n \le p_{\text{I}} + p_{\text{II}}.\]
In region $\text{II}$, note that the hypothesis of Lemma \ref{new-blocks} is satisfied since if $B_{l,j}$ and $B_{l,j-1}$ are two consecutive blocks in that region $|B_{l-1,j-1}|\ge(\sqrt{2n}+1)\sqrt{2n}/2$. By Lemma \ref{new-blocks}, a subblock $B$ of $B_{l,j}B_{l,j-1}$ is a new $n$-block at level $l+1$ if and only if $a^{l-j}b^{j-1}$ is a subblock of $B$. Therefore, the number of new $n$-subblocks of $B_{l,j}B_{l,j-1}$ is equal to $n-l$. Thus \begin{align*}p_\text{II}&=\sum_{\sqrt{2n}<l<n+1}(n-l)(l-3)\\
&=n\,\sum_{l=1}^{n+1}l - \sum_{l=1}^{n+1}l^2 + o(n^3)\\
&=n\,\fr{(n+1)(n+2)}{2}-\fr{(n+1)(n+2)(2n+3)}{6}+o(n^3)\\
&=\fr{n^3}{2}-\fr{2n^3}{6}+o(n^3)=\fr{n^3}{6}+o(n^3).\end{align*}
For region I, a coarse approximation gives \[p_\text{I}\le \sum_{1\le l\le \sqrt{2n}}(n-1)(l-1) + 4\sum_{\sqrt{2n}<l\le n+1}(l-1)=o(n^3).\] It follows that \[\frac{n^3}{6}+o(n^3)\le p_n \le \frac{n^3}{6}+o(n^3).\]

\begin{figure}[hbt]
\begin{picture}(0,0)%
\includegraphics{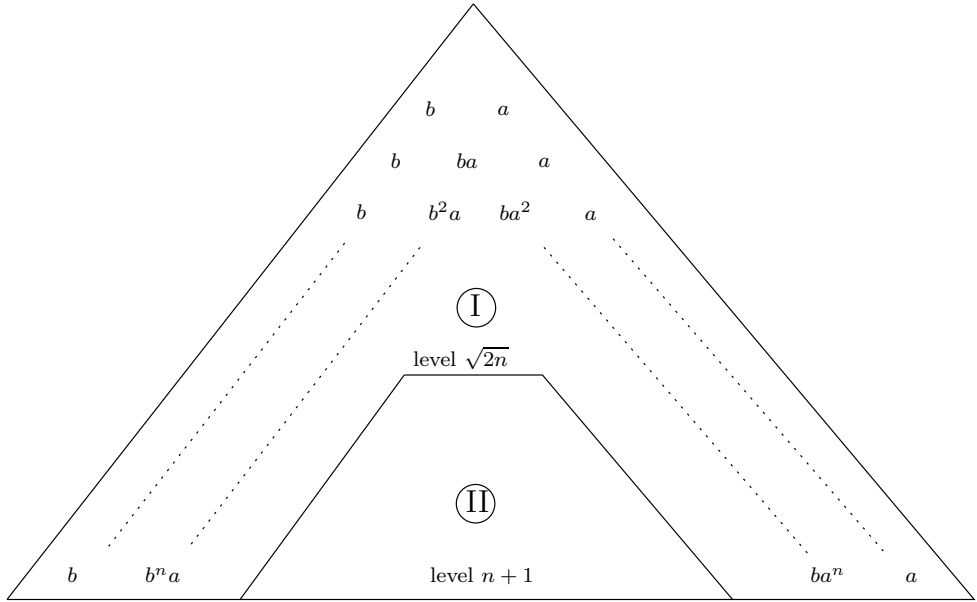}%
\end{picture}%
\setlength{\unitlength}{3947sp}%
\begingroup\makeatletter\ifx\SetFigFont\undefined%
\gdef\SetFigFont#1#2#3#4#5{%
  \reset@font\fontsize{#1}{#2pt}%
  \fontfamily{#3}\fontseries{#4}\fontshape{#5}%
  \selectfont}%
\fi\endgroup%
\begin{picture}(6099,3768)(2389,-2542)
\put(5222,183){\makebox(0,0)[lb]{\smash{\SetFigFont{8}{9.6}{\familydefault}{\mddefault}{\updefault}{\color[rgb]{0,0,0}$ba$}%
}}}
\put(2781,-2422){\makebox(0,0)[lb]{\smash{\SetFigFont{8}{9.6}{\familydefault}{\mddefault}{\updefault}{\color[rgb]{0,0,0}$b$}%
}}}
\put(8042,-2422){\makebox(0,0)[lb]{\smash{\SetFigFont{8}{9.6}{\familydefault}{\mddefault}{\updefault}{\color[rgb]{0,0,0}$a$}%
}}}
\put(4951,-1065){\makebox(0,0)[lb]{\smash{\SetFigFont{8}{9.6}{\familydefault}{\mddefault}{\updefault}{\color[rgb]{0,0,0}level $\sqrt{2n}$}%
}}}
\put(3269,-2422){\makebox(0,0)[lb]{\smash{\SetFigFont{8}{9.6}{\familydefault}{\mddefault}{\updefault}{\color[rgb]{0,0,0}$b^{n}a$}%
}}}
\put(5059,-2422){\makebox(0,0)[lb]{\smash{\SetFigFont{8}{9.6}{\familydefault}{\mddefault}{\updefault}{\color[rgb]{0,0,0}level $n+1$}%
}}}
\put(7446,-2422){\makebox(0,0)[lb]{\smash{\SetFigFont{8}{9.6}{\familydefault}{\mddefault}{\updefault}{\color[rgb]{0,0,0}$ba^{n}$}%
}}}
\put(5493,-142){\makebox(0,0)[lb]{\smash{\SetFigFont{8}{9.6}{\familydefault}{\mddefault}{\updefault}{\color[rgb]{0,0,0}$ba^2$}%
}}}
\put(4595,-136){\makebox(0,0)[lb]{\smash{\SetFigFont{8}{9.6}{\familydefault}{\mddefault}{\updefault}{\color[rgb]{0,0,0}$b$}%
}}}
\put(4811,183){\makebox(0,0)[lb]{\smash{\SetFigFont{8}{9.6}{\familydefault}{\mddefault}{\updefault}{\color[rgb]{0,0,0}$b$}%
}}}
\put(5028,508){\makebox(0,0)[lb]{\smash{\SetFigFont{8}{9.6}{\familydefault}{\mddefault}{\updefault}{\color[rgb]{0,0,0}$b$}%
}}}
\put(5480,514){\makebox(0,0)[lb]{\smash{\SetFigFont{8}{9.6}{\rmdefault}{\mddefault}{\updefault}{\color[rgb]{0,0,0}$a$}%
}}}
\put(6028,-142){\makebox(0,0)[lb]{\smash{\SetFigFont{8}{9.6}{\familydefault}{\mddefault}{\updefault}{\color[rgb]{0,0,0}$a$}%
}}}
\put(5738,189){\makebox(0,0)[lb]{\smash{\SetFigFont{8}{9.6}{\familydefault}{\mddefault}{\updefault}{\color[rgb]{0,0,0}$a$}%
}}}
\put(5047,-142){\makebox(0,0)[lb]{\smash{\SetFigFont{8}{9.6}{\familydefault}{\mddefault}{\updefault}{\color[rgb]{0,0,0}$b^2a$}%
}}}
\put(5306,-746){\makebox(0,0)[lb]{\smash{\SetFigFont{12}{14.4}{\familydefault}{\mddefault}{\updefault}{\color[rgb]{0,0,0}I}%
}}}
\put(5282,-1981){\makebox(0,0)[lb]{\smash{\SetFigFont{12}{14.4}{\familydefault}{\mddefault}{\updefault}{\color[rgb]{0,0,0}II}%
}}}
\end{picture}
\caption{In region II, two consecutive basic blocks $B_{l,j}$ and $B_{l,j-1}$ at level $l$ will generate $n-l$ new $n$-blocks at level $l+1$, according to Lemma \ref{new-blocks}.}
\label{fig:block-count}
\end{figure}
\end{proof}

An immediate corollary is that the countable-substitution subshift $(\Sigma,\sigma)$ has topological entropy \[h_\text{top}(\Sigma,\sigma)=\lim_{n\to\infty}\fr{\log(p_n)}{n}=0.\] The variational principle implies that the measure-theoretic entropy $h(\sigma)$ with respect to every measure $\nu_\alpha$ is zero. Since entropy is an isomorphism invariant, this proves that the Pascal adic has zero entropy for every Bernoulli measure $\mu_\alpha$.



\subsection{Directional unique ergodicity}

For any two blocks $B$ and $C$ in the language $\mathcal{L}(\zeta)$, let $\freq(B,C)$ denote the frequency of occurrences of $B$ in $C$. In particular, if $\omega\in\Sigma$, and if $\omega_{-N}^N$ denotes the block $\omega_{-N}\omega_{-N+1}\dots\omega_0\dots\omega_{N-1}\omega_N$, then, provided that $2N+1$ is greater than $|B|$ (the length of $B$), we have  \[\freq(B,\omega_{-N}^N)=\fr{1}{2N+1}\sum_{i=-N}^{N-|B|}\II_{[B]}\circ\sigma^i(\omega).\]
The ergodicity of each $\nu_{\alpha}=\phi(\mu_{\alpha})$ (carried by the isomorphism) implies therefore that for every block $B\in\mathcal{L}(\zeta)$ and $\nu_{\alpha}$-a.e. $\omega\in\Sigma$
\begin{equation}\label{equ:freq1}
\lim_{N\to\infty}\freq(B,\omega_{-N}^N)=\nu_{\alpha}(B).
\end{equation}
Let $x$ be a path in the Pascal graph going through the vertices $(n,k_n(x))$, and let $\omega=\phi(x)$. If $x$ is not eventually diagonal, then there are sequences $i_n(x)\nearrow\infty$ and $j_n(x)\nearrow\infty$ such that $B_{n,k_n(x)}=\omega_{-i_n(x)}^{j_n(x)}$, where $B_{n,k_n}$ is the basic block at vertex $(n,k_n)$ in the Pascal triangle of words. Therefore, (\ref{equ:freq1}) implies that for $\mu_\alpha$-a.e. $x$ \[\lim_{n\to\infty}\freq(B,B_{n,k_n(x)})=\nu_{\alpha}(B).\]
In other words, when $(n,k_n)$ for $n=1,2,\dots$ are the vertices of a generic path, which goes down the Pascal graph at an ``angle'' $\alpha$,
we have $\freq(B,B_{n,k_n}) \to \nu_\alpha (B)$.
We strengthen this statement as follows:

\begin{thm}\label{thm:freq} For any block $B\in\mathcal{L}(\zeta)$,
any $\alpha\in(0,1)$, and \emph{any} sequence $k_n\to \infty$ such that $k_n/n\to\alpha$ we have \[\lim_{n\to\infty}\freq(B,B_{n,k_n})=\nu_{\alpha}(B).\]\label{equ:freq2}
\end{thm}
\begin{figure}[hbt]
\begin{picture}(0,0)%
\includegraphics{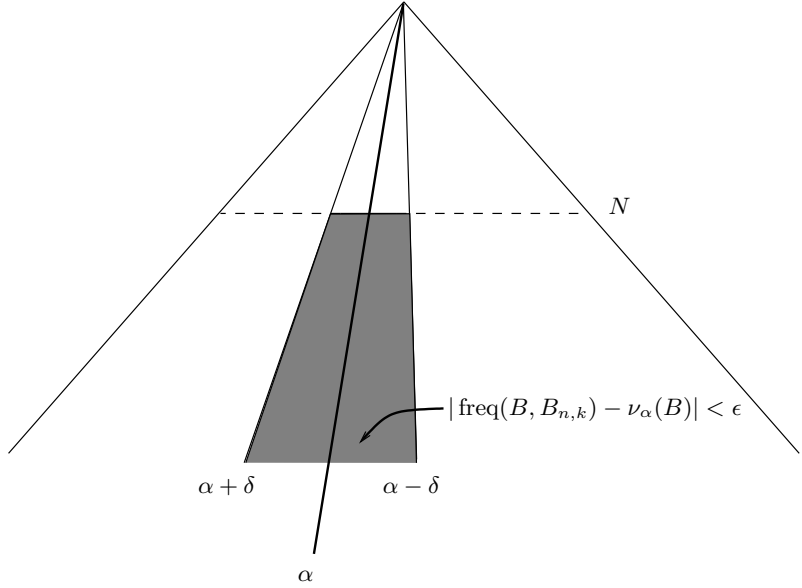}%
\end{picture}%
\setlength{\unitlength}{3947sp}%
\begingroup\makeatletter\ifx\SetFigFont\undefined%
\gdef\SetFigFont#1#2#3#4#5{%
  \reset@font\fontsize{#1}{#2pt}%
  \fontfamily{#3}\fontseries{#4}\fontshape{#5}%
  \selectfont}%
\fi\endgroup%
\begin{picture}(4990,3700)(2689,-4762)
\put(6470,-2404){\makebox(0,0)[lb]{\smash{\SetFigFont{9}{10.8}{\rmdefault}{\bfdefault}{\updefault}{\color[rgb]{0,0,0}$N$}%
}}}
\put(3891,-4161){\makebox(0,0)[lb]{\smash{\SetFigFont{9}{10.8}{\rmdefault}{\bfdefault}{\updefault}{\color[rgb]{0,0,0}$\alpha+\delta$}%
}}}
\put(5051,-4161){\makebox(0,0)[lb]{\smash{\SetFigFont{9}{10.8}{\rmdefault}{\bfdefault}{\updefault}{\color[rgb]{0,0,0}$\alpha-\delta$}%
}}}
\put(4517,-4717){\makebox(0,0)[lb]{\smash{\SetFigFont{9}{10.8}{\rmdefault}{\bfdefault}{\updefault}{\color[rgb]{0,0,0}$\alpha$}%
}}}
\put(5459,-3669){\makebox(0,0)[lb]{\smash{\SetFigFont{9}{10.8}{\rmdefault}{\bfdefault}{\updefault}{\color[rgb]{0,0,0}$|\freq(B,B_{n,k})-\nu_\alpha(B)|<\epsilon$}%
}}}
\end{picture}
\caption{For any block $B_{n,k}$ in the filled area the frequency of appearances of $B$ in $B_{n,k}$ is $\epsilon$-close to $\nu_{\alpha}(B)$ (for $\delta$ small and $N$ large).}
\label{fig:directional_ergodicity}
\end{figure}

To prove Theorem \ref{thm:freq} we determine explicitly where a given block $B$ is made and how many times it appears in each basic block. First we introduce some notation and recall the key structure of the basic blocks. For a fixed  block $B\in\mathcal{L}(\zeta)$, denote by $a(B,n,k)$ the number of occurrences of $B$ in $B_{n,k}$. If no confusion is possible we will simply denote it by $a(n,k)$, so that \begin{equation}\label{eq:freq}\freq(B,B_{n,k_n})=\fr{a(n,k_n)}{|B_{n,k_n}|}=\fr{a(n,k_n)}{C(n,k_n)}.\end{equation}
By induction it is easy to show that \[B_{n,k}=\underbrace{\boldsymbol{b^k}a\dots b\boldsymbol{a^{n-k-1}}}_{B_{n-1,k}}\underbrace{\boldsymbol{b^{k-1}}a\dots b\boldsymbol{a^{n-k}}}_{B_{n-1,k-1}}.\] This decomposition characterizes the basic blocks $B_{n,k}$; if one sees the \emph{telltale} block $b{\boldsymbol{a^{n-k-1}b^{k-1}}}a$, it is always within $B_{n,k}$, at the ``join'', as seen above. The previous structure is easily seen by writing only the beginning and ending of the basic blocks in the Pascal triangle of words; for example, at level $n$ we have: \[\underbrace{\boldsymbol{b}}_{B_{n,n}}\quad \underbrace{\boldsymbol{b^{n-1}a}}_{B_{n,n-1}} \quad \dots \quad \underbrace{{\boldsymbol{b^k}}a\dots b{\boldsymbol{a^{n-k}}}}_{B_{n,k}} \quad \dots \quad \underbrace{{\boldsymbol{ba^{n-1}}}}_{B_{n,1}} \quad\underbrace{{\boldsymbol{a}}}_{B_{n,0}}\] Figure \ref{fig:structure_diagram} shows how the structure is carried from one level to the next:\begin{figure}[ht]\begin{gather*} \overbrace{{\boldsymbol{b^{k-1}}}a\dots b{\boldsymbol{ a^{n-k-1}}}}^{B_{n-2,k-1}}\\ \underbrace{{\boldsymbol{b^{k}}}a\dots b{\boldsymbol{ a^{n-k-1}}}}_{B_{n-1,k}} \quad \underbrace{{\boldsymbol{b^{k-1}}}a\dots b{\boldsymbol{ a^{n-k}}}}_{B_{n-1,k-1}} \\ {\boldsymbol{b^{k+1}}}a\dots b{\boldsymbol{ a^{n-k-1}}} \quad \underbrace{{\boldsymbol{b^{k}}}a\dots b{\boldsymbol{ a^{n-k}}}}_{B_{n,k}} \quad {\boldsymbol{b^{k-1}}}a\dots b{\boldsymbol{ a^{n-k+1}}} \\{\boldsymbol{b^{k+2}}}a\dots b{\boldsymbol{ a^{n-k-1}}} \quad  {\boldsymbol{b^{k+1}}}a\dots b{\boldsymbol{ a^{n-k}}} \quad {\boldsymbol{b^{k}}}a\dots b{\boldsymbol{ a^{n-k+1}}} \quad {\boldsymbol{b^{k-1}}}a\dots b{\boldsymbol{ a^{n-k+2}}}\end{gather*}\caption{}\label{fig:structure_diagram}\end{figure}

\begin{lemma}\label{lem:first_appearance}
Let $B\in\mathcal{L}(\zeta)$. Then there is a unique vertex $(n_0,k_0)$ such that $B$ is a subblock of $B_{n_0,k_0}$ and $B$ does not appear in any other basic block in the ``rectangle'' above $(n_0,k_0)$ formed by all vertices $(n,k)$ with $k<k_0$ and $n-k<n_0-k_0$.
\end{lemma}

\begin{proof}
Assume that $B$ is a subblock of both $B_{n_0,k_0}$ and $B_{n_0',k_0'}$, and that $B$ does not appear in any other basic block in the ``rectangles'' above $(n_0,k_0)$ and $(n_0',k_0')$ (in particular we can assume that $n_0'>n_0$ and $k_0'>k_0$) --- see Figure \ref{fig:first_appearance2}.  Recall that \begin{align*}& B_{n_0,k_0}=b^{k_0} a\dots ba^{n_0-k_0-1}|b^{k_0-1}\dots ba^{n_0-k_0}\\ & B_{n_0',k_0'}=b^{k_0'}a\dots ba^{n_0'-k_0'-1}|b^{j_0'-1}\dots ba^{n_0'-k_0'},\end{align*} so that $B$ must contain a central subblock of both \begin{align*} & \quad\quad \, b\overbrace{a\dots a}^{n_0'-k_0'-1}|\overbrace{b \dots\dots\dots b}^{k_0'-1}a \\ & b\underbrace{a\dots\dots\dots a\;}_{n_0-k_0-1}|\underbrace{b \dots b}_{k_0-1}a. \end{align*} (By central we mean containing $a|b$, where ``$|$'' symbolizes the splitting into blocks of the above level). The only way it can happen is if $B=a^ib^j$ with $i<n_0'-k_0'-1$ and $j<k_0-1$.  Therefore $B$ had to appear at the vertex $(n_0'-k_0'+k_0,k_0)$ whose basic block is $\dots ba^{n_0'-k_0'-1}b^{k_0-1}a\dots$, contradicting our initial assumption.
\begin{figure}[hbt]
\begin{picture}(0,0)%
\includegraphics{first_appearance2.pstex}%
\end{picture}%
\setlength{\unitlength}{3947sp}%
\begingroup\makeatletter\ifx\SetFigFont\undefined%
\gdef\SetFigFont#1#2#3#4#5{%
  \reset@font\fontsize{#1}{#2pt}%
  \fontfamily{#3}\fontseries{#4}\fontshape{#5}%
  \selectfont}%
\fi\endgroup%
\begin{picture}(6945,3287)(1993,-3636)
\put(5868,-2399){\makebox(0,0)[lb]{\smash{\SetFigFont{7}{8.4}{\familydefault}{\mddefault}{\updefault}{\color[rgb]{0,0,0}$(n_0,k_0)$}%
}}}
\put(3931,-3520){\makebox(0,0)[lb]{\smash{\SetFigFont{7}{8.4}{\familydefault}{\mddefault}{\updefault}{\color[rgb]{0,0,0}$(n_0',k_0')$}%
}}}
\put(5611,-3113){\makebox(0,0)[lb]{\smash{\SetFigFont{7}{8.4}{\familydefault}{\mddefault}{\updefault}{\color[rgb]{0,0,0}$(n_0'-k_0'+k_0,k_0)$}%
}}}
\put(2326,-661){\makebox(0,0)[lb]{\smash{\SetFigFont{7}{8.4}{\familydefault}{\mddefault}{\updefault}{\color[rgb]{0,0,0}Region where we assume $B$ does not appear.}%
}}}
\end{picture}
\label{fig:first_appearance2}
\end{figure}
\end{proof}

\begin{remark} From the uniqueness in Lemma \ref{lem:first_appearance}, it follows that if the block $B$ appears for the first time at the vertex $(n_0,k_0)$, then it will only appear as a subblock in the triangle below the vertex $(n_0,k_0)$, as shown in Figure \ref{fig:first_appearance}.
\begin{figure}[hbt]
\begin{picture}(0,0)%
\includegraphics{first_appearance.pstex}%
\end{picture}%
\setlength{\unitlength}{3947sp}%
\begingroup\makeatletter\ifx\SetFigFont\undefined%
\gdef\SetFigFont#1#2#3#4#5{%
  \reset@font\fontsize{#1}{#2pt}%
  \fontfamily{#3}\fontseries{#4}\fontshape{#5}%
  \selectfont}%
\fi\endgroup%
\begin{picture}(6624,3174)(2089,-3523)
\put(5476,-2011){\makebox(0,0)[lb]{\smash{\SetFigFont{8}{9.6}{\familydefault}{\mddefault}{\updefault}{\color[rgb]{0,0,0}$(n_0,k_0)$}%
}}}
\put(7431,-1662){\makebox(0,0)[lb]{\smash{\SetFigFont{12}{14.4}{\familydefault}{\mddefault}{\updefault}{\color[rgb]{0,0,0}All vertices $(n_0+i,k_0+j)$ }%
}}}
\put(7579,-1990){\makebox(0,0)[lb]{\smash{\SetFigFont{12}{14.4}{\familydefault}{\mddefault}{\updefault}{\color[rgb]{0,0,0}for $i\ge 0$ and $0\le j \le i$.}%
}}}
\put(2476,-661){\makebox(0,0)[lb]{\smash{\SetFigFont{8}{9.6}{\familydefault}{\mddefault}{\updefault}{\color[rgb]{0,0,0}Region where the block $B$ appears.}%
}}}
\end{picture}
\label{fig:first_appearance}
\end{figure}
\end{remark}

\begin{lemma}\label{lem:atmosttwice}
If a block $B$ appears in $B_{n_0,k_0}$ for the first time, then it can appear at most twice (i.e $a(n_0,k_0)=1$ or $2$).

\end{lemma}

\begin{proof}
If $B$ appears for the first time in $B_{n_0,k_0}$, then it is a central subblock of \[b^{k_0}\dots\underbrace{\dots a^{n_0-k_0-1}|b^{k_0-1}\dots}_{B}\dots a^{n_0-k_0}.\] There are two possibilities: either $B$ contains the telltale block $ba^{n_0-k_0-1}|b^{k_0-1}a,$ or it does not. In the first case, $B$ cannot appear twice, since the telltale block appears only once in $B_{n_0,k_0}$. In the second case, the only way that $B$ could appear twice is if it would start with $a^ib^{k_0-1}$ for some $i\le n_0-k_0-1$, and end with $a^{n_0-k_0-1}b^j$ for some $j\le k_0-1$: \begin{align*}& b^{k_0}\dots\dots \overbrace{a^{n_0-k_0-1}|b^{k_0-1}\dots}^{B}\dots a^{n_0-k_0} \\ & b^{k_0}\dots\underbrace{\dots a^{n_0-k_0-1}|b^{k_0-1}}_{B}\dots\dots a^{n_0-k_0}.\end{align*} Then clearly $B$ cannot appear a third time.
\end{proof}

Given a block $B\in\mathcal{L}(\zeta)$, assume that $B$ appears for the first time at vertex $(n_0,k_0)$ in the Pascal triangle of words (this first appearance is unique by Lemma \ref{lem:first_appearance}). Define the \emph{triangle of appearances of $B$} to be \begin{gather*} a(n_0,k_0) \\ a(n_0+1,k_0+1)\quad a(n_0+1,k_0) \\ a(n_0+2,k_0+2)\quad a(n_0+2,k_0+1)\quad a(n_0+2,k_0) \\ \dots \quad \dots \quad \dots \quad \dots \quad \dots \quad \dots \quad \dots \quad \dots \quad \dots \quad \dots\quad \dots\end{gather*}
By Lemma \ref{lem:first_appearance}, all other frequencies are zero, namely $a(n,k)=0$ for all vertices $(n,k)$ not equal to $(n_0+i,k_0+j)$ where $i\ge 0$ and $0\le j\le i$ --- see Figure \ref{fig:first_appearance}.



\begin{example}\label{ex:triangle_of_appearances}
The triangle of appearances of $B=a^2b^2$ is \begin{gather*} 1 \\ 2 \quad 2 \\ 3 \quad 5 \quad 2 \\ 4 \quad 9 \quad 9 \quad 4 \\  \dots \quad \dots \quad \dots \end{gather*} where the first element in the triangle is $a(6,3)$.
\end{example}


It turns out that the triangle of appearances of a block $B$ is either of the type of example \ref{ex:triangle_of_appearances} or the Pascal triangle with certain ``initial diagonals'' removed.

\begin{lemma}\label{lem:fivecases}
Let $B\in\mathcal{L}(\zeta)$ and assume $B$ appears for the first time in $B_{n_0,k_0}$. Let $i_0=n_0-k_0-1$ and $j=k_0-1$. Then there are exactly five possible cases for the triangle of appearances of $B$: \\
\underline{Case 1}: $B=\boldsymbol{a^{i_0}b^{j_0}}$. The number of occurrences of $B$ satisfies the recurrence relation \[a(n,k)=a(n-1,k)+a(n-1,k-1)+1,\] and the triangle of occurrences of $B$ is \begin{gather*} 1 \\ 2 \quad 2 \\ 3 \quad 5 \quad 2 \\ 4 \quad 9 \quad 9 \quad 4 \\  \dots \quad \dots \quad \dots \end{gather*} It follows that $a(n,k)=C(n-n_0+2,k-k_0+1)-1.$\\
\underline{Case 2}: $B={\dots}b\boldsymbol{a^{i_0}b^{j_0}}a{\dots}$ ($B$ contains the telltale block). The number of occurrences of $B$ satisfies the recurrence relation \[a(n,k)=a(n-1,k)+a(n-1,k-1),\] and the triangle of occurrences of $B$ is just the Pascal triangle. \begin{gather*} 1 \\ 1 \quad 1 \\ 1 \quad 2 \quad 1 \\ 1 \quad 3 \quad 3 \quad 1 \\  \dots \quad \dots \quad \dots \end{gather*} It follows that $a(n,k)=C(n-n_0,k-k_0)$. \\
\underline{Case 3}: $B=\boldsymbol{a^iB_{n_0-2,k_0-1}b^j}$ for some $i<i_0$ and $j<j_0$. The number of occurrences of $B$ satisfies the recurrence relation \[a(n,k)=a(n-1,k)+a(n-1,k-1),\] and the triangle of occurrences of $B$ is \begin{gather*} 2 \\ 3 \quad 3 \\ 4 \quad 6 \quad 4 \\ 5 \quad 10 \quad 10 \quad 5 \\  \dots \quad \dots \quad \dots \end{gather*} It follows that $a(n,k)=C(n-n_0+2,k-k_0+1)$.\\
\underline{Case 4}: $B={\dots}b\boldsymbol{a^{i_0}b^j}$ for some $j\le j_0$, and $B$ is not as in Case 3. The number of occurrences of $B$ satisfies the recurrence relation \[a(n,k)=a(n-1,k)+a(n-1,k-1),\] and the triangle of occurrences of $B$ is the following \begin{gather*} 1 \\ 2 \quad 1 \\ 3 \quad 3 \quad 1 \\ 4 \quad 6 \quad 4 \quad 1 \\  \dots \quad \dots \quad \dots \end{gather*} It follows that $a(n,k)=C(n-n_0+1,k-k_0)$.\\
\underline{Case 5}: $B=\boldsymbol{a^{i}b^{j_0}}a{\dots}$ for some $i\le i_0$, and $B$ is not as in Case 3. The number of occurrences of $B$ satisfies the recurrence relation \[a(n,k)=a(n-1,k)+a(n-1,k-1),\] and the triangle of occurrences of $B$ is the following \begin{gather*} 1 \\ 1 \quad 2 \\ 1 \quad 3 \quad 3 \\ 1 \quad 4 \quad 6 \quad 4 \\  \dots \quad \dots \quad \dots \end{gather*} It follows that $a(n,k)=C(n-n_0+1,k-k_0+1)$.
\end{lemma}

\begin{proof}
If $B$ appears for the first time in $B_{n_0,k_0}$, then clearly $B$ is one of the following: \ben \item $B=\boldsymbol{a^{i_0}b^{j_0}}$ \item $B={\dots}b\boldsymbol{a^{i_0}b^{j_0}}a{\dots}$ \item $B={\dots}b\boldsymbol{a^{i_0}b^j}$ for some $j\le j_0$ \item $B=\boldsymbol{a^{i}b^{j_0}}a{\dots}$ for some $i\le i_0$. \een

In Case (i), which corresponds to Case 1, the block $B$ appears only once in $B_{n_0,k_0}$, and is made infinitely many times further down. More precisely, every concatenation of two adjacent basic blocks below $n_0$ creates a unique block $B$. It follows that $a(n,k)=a(n-1,k)+a(n-1,k-1)+1$, and the triangle of occurrences of $B$ is \begin{gather*} 1 \\ 2 \quad 2 \\ 3 \quad 5 \quad 2 \\ 4 \quad 9 \quad 9 \quad 4 \\  \dots \quad \dots \quad \dots \end{gather*} In the previous triangle, $a(n+n_0,k+k_0)$ represents the number of finite paths from the root $(0,0)$ to the vertex $(n,k)$ in a modified Pascal graph with extra ``wormholes'' --- paths connecting directly $(0,0)$ to each $(n,k)$  --- see Figure \ref{fig:wormhole}. It is not hard to see from this graph that $a(n+n_0,k+k_0)+1$ is equal to the sum of the binomial coefficients inside the rectangle determined by $(0,0)$ and $(n,k)$, i.e. \[a(n+n_0,k+k_0)+1= \sum_{i=0}^{n-k}\sum_{j=i}^{k+i}C(j,j-i).\] Using well-known properties of the binomial coefficients we get \begin{align*} a(n+n_0,k+k_0)+1&= \sum_{i=0}^{n-k}\sum_{j=i}^{k+i}C(j,j-i)\\ &= \sum_{i=0}^{n-k} C(k+i+1,k+1)\\ &=C(n+2,k+1),\end{align*} showing as announced that $a(n,k)=C(n-n_0+2,k-k_0+1)-1.$
\begin{figure}[hbt]
\begin{picture}(0,0)%
\includegraphics{wormhole.pstex}%
\end{picture}%
\setlength{\unitlength}{3947sp}%
\begingroup\makeatletter\ifx\SetFigFont\undefined%
\gdef\SetFigFont#1#2#3#4#5{%
  \reset@font\fontsize{#1}{#2pt}%
  \fontfamily{#3}\fontseries{#4}\fontshape{#5}%
  \selectfont}%
\fi\endgroup%
\begin{picture}(5658,3038)(2755,-3634)
\put(5433,-692){\makebox(0,0)[lb]{\smash{\SetFigFont{7}{8.4}{\rmdefault}{\bfdefault}{\updefault}{\color[rgb]{0,0,0}$(0,0)$}%
}}}
\end{picture}
\label{fig:wormhole}
\end{figure}

In Case (ii), a similar argument as in the proof of Lemma \ref{lem:atmosttwice} shows that  $B$ can be made only once. Therefore $a(n,k)=a(n-1,k)+a(n-1,k-1)$, and the triangle of appearances of $B$ is just the Pascal triangle. Since the first element in the triangle is $a(n_0,k_0)$, it follows that $a(n,k)=C(n-n_0,k-k_0)$, establishing Case 2.

A special case of (iii) and (iv) is Case 3, i.e when $B=\boldsymbol{a^iB_{n_0-2,k_0-1}b^j}$ for some $i<i_0$ and $j<j_0$. Here, $B$ is made twice at $(n_0,k_0)$, and infinitely many times along the edges of the triangle below $(n_0,k_0)$. First, observe that \[B_{n_0,k_0}=b^{k_0}\dots a^{n_0-k_0-2}B_{n_0-2,k_0-1}|B_{n_0-2,k_0-1}b^{k_0-2}\dots \dots a^{n_0-k_0}.\] Therefore $B$ appears twice in $B_{n_0,k_0}$: \begin{align*} &b^{k_0}\dots \overbrace{a^{n_0-k_0-2}B_{n_0-2,k_0-1}|b^{k_0-1}}^{B}\dots  \dots a^{n_0-k_0} \\ &b^{k_0}\dots\dots \underbrace{a^{n_0-k_0-1}|B_{n_0-2,k_0-1}b^{k_0-2}}_{B}\dots \dots a^{n_0-k_0}. \end{align*} Second, note that $B$ can be made exactly once at each vertex along the edges of the triangle below $(n_0,k_0)$: \[B_{n_0+j,k_0}=\dots \overbrace{a^{n_0-k_0+j}|B_{n_0-2,k_0-1}b^{k_0-2}}^{B}\dots,\] and \[ B_{n_0+j,k_0+j}=\dots \underbrace{a^{n_0-k_0-2}B_{n_0-2,k_0-1}|b^{k+j}}_{B}\dots,\] for all $j\ge 1$. Thus $a(n,k)=a(n-1,k)+a(n-1,k-1)$, and the triangle of occurrences of $B$ is \begin{gather*} 2 \\ 3 \quad 3 \\ 4 \quad 6 \quad 4 \\ 5 \quad 10 \quad 10 \quad 5 \\  \dots \quad \dots \quad \dots \end{gather*} which is a subtriangle of the Pascal triangle, hence the relation $a(n,k)=C(n-n_0+1,k-k_0)$ follows.

To conclude, assume that $B$ is as in (iii) but not of the type $B=a^iB_{n_0-2,k_0-1}b^j$. Then $B=\dots b a^{n_0-k_0-1}b^j$ appears only once in $B_{n_0,k_0}$, because in order to appear twice it would have to contain the entire block $B_{n_0-2,k_0-1}$. For the same reason it cannot be made again along the vertices $(n_0+j,k_0)$ for $j\ge 1$. Any two adjacent basic blocks inside the triangle below $(n_0,k_0)$ will not produce any more $B$'s either, because they are all of the type $\dots ba^{i}|b^{j}a\dots$ where $i>i_0$. On the other hand,  $B$ is created once more at each vertex $(n_0+j,k_0+j)$ for $j\ge 1$: \[B_{n_0+j,k_0+j}=\dots\underbrace{\dots ba^{n_0-k_0-1}|b^{k_0+j}}_{B}\dots.\] Therefore the triangle of appearances is \begin{gather*} 1 \\ 2 \quad 1 \\ 3 \quad 3 \quad 1 \\ 4 \quad 6 \quad 4 \quad 1 \\  \dots \quad \dots \quad \dots \end{gather*} As previously this is a subtriangle of the Pascal triangle, thus $a(n,k)=C(n-n_0+1,k-k_0)$, and Case 4 is established. Case 5 is similar.
\end{proof}

\begin{lemma}\label{lem:quotientbino}
Let $r\ge 0$ and $0\le s \le r$. If $k_n/n\to\alpha$ as $n\to\infty$, then \[\lim_{n\to\infty}\fr{C(n-r,k_n-s)}{C(n,k_n)}=\alpha^s(1-\alpha)^{r-s}.\]
\end{lemma}


\begin{remark}
This was the key calculation used by Hajian, Ito, and Kakutani to prove the ergodicity of $\mu_\alpha$ for the Pascal adic --- see \cite{hik}, and also \cite[Theorem 2.7]{gibbs}.
\end{remark}

\begin{proof}[Proof of Theorem \ref{thm:freq}]
By (\ref{eq:freq}) we need to show that \begin{equation}\label{eq:freq2}\lim_{n\to\infty}\fr{a(n,k_n)}{C(n,k_n)}=\nu_{\alpha}(B),\end{equation} when $k_n/n\to\alpha$. According to Lemma \ref{lem:fivecases}, $a(n,k_n)$ is equal to one of the following: \ben \item[Case 1:] $C(n-n_0+2,k-k_0+1)-1$ \item[Case 2:] $C(n-n_0,k-k_0)$ \item[Case 3:] $C(n-n_0+2,k-k_0+1)$ \item[Case 4:] $C(n-n_0+1,k-k_0)$ \item[Case 5:] $C(n-n_0+1,k-k_0+1)$.\een In each of these cases Lemma \ref{lem:quotientbino} shows that the limit in (\ref{eq:freq2}) equals respectively: \ben \item[Case 1:] $\alpha^{k_0-1}(1-\alpha)^{n_0-k_0-1}$ \item[Case 2:] $\alpha^{k_0}(1-\alpha)^{n_0-k_0}$ \item[Case 3:] $\alpha^{k_0-1}(1-\alpha)^{n_0-k_0-1}$ \item[Case 4:] $\alpha^{k_0}(1-\alpha)^{n_0-k_0-1}$ \item[Case 5:] $\alpha^{k_0-1}(1-\alpha)^{n_0-k_0}$.\een

On the other hand, $\nu_\alpha(B)=\sum\nu_\alpha(B_{n,k}),$ where the sum is taken over all $(n,k)$ such that $B$ is created in the basic block $B_{n,k}$, and if $B$ is created twice in $B_{n,k}$, then $\nu_\alpha(B_{n,k})$ appears twice in the sum. In particular, if $B=B_{n,k}$ we have $\nu_\alpha(B)=\alpha^{k}(1-\alpha)^{n-k}$. So for example, if $B$ is as in Case 3, we have:
\begin{align*}\nu_\alpha(B)&=2\nu_\alpha(B_{n_0,k_0})+\sum_{i=1}^{\infty}\nu_\alpha(B_{n_0+i,k_0+i})+\sum_{i=1}^{\infty}\nu_\alpha(B_{n_0+i,k_0})\\&=\sum_{i=0}^{\infty}\nu_\alpha(B_{n_0+i,k_0+i})+\sum_{i=0}^{\infty}\nu_\alpha(B_{n_0+i,k_0})\\&= \sum_{i=0}^{\infty}\alpha^{k_0+i}(1-\alpha)^{n_0-k_0} + \sum_{i=0}^{\infty}\alpha^{k_0}(1-\alpha)^{n_0-k_0+i}\\&=\alpha^{k_0}(1-\alpha)^{n_0-k_0}\left[\sum_{i=0}^{\infty} \alpha^i + \sum_{i=0}^{\infty} (1-\alpha)^i \right] \\ &= \alpha^{k_0}(1-\alpha)^{n_0-k_0}\left[\fr{1}{1-\alpha}+\fr{1}{\alpha}\right] \\&= \alpha^{k_0}(1-\alpha)^{n_0-k_0}\;\fr{\alpha + 1-\alpha}{\alpha(1-\alpha)}\\&=\alpha^{k_0-1}(1-\alpha)^{n_0-k_0-1}. \end{align*}
Similar computations show that the identity (\ref{eq:freq2}) also holds in the other cases.
\end{proof}

\begin{remark}
We could have used the Ergodic Theorem to show that the limit in
(\ref{eq:freq2}), which is known to exist because of Lemma
\ref{lem:fivecases}, must be equal to $\nu_\alpha(B)$: as noted
previously (p. \pageref{equ:freq2}), for $\mu_\alpha$-a.e path $x$
we have $k_n(x)/n\to\alpha$ and
\[\lim_{n\to\infty}\fr{a(n,k_n(x))}{C(n,k_n)}=\nu_\alpha(B).\]
\end{remark}

\begin{remark}
We thank Anthony Quas for an insight which leads to a proof simpler
than the preceding one.
It is possible to show that that after a while (that is, below
a certain level) ``new'' creations of a block $B$ due to
concatenations constitute a negligible fraction of the total number of
appearances of $B$.  More specifically, fix an allowed block $B\in
\mathcal L (\Sigma)$ and a positive
integer $m$ large enough that $B$ appears at level $m$ of the Pascal
triangle of words.
Consider a level $n$ with
$n\geq m$. As mentioned above, any basic block $B_{n,k}$ at this level
factors uniquely into basic blocks of level $m$. Let us call ``old''
appearances of $B$ in $B_{n,k}$ the ones that are contained entirely in
one of these basic blocks of level $m$, and ``new'' appearances the
rest. The new appearances were created by concatenations after level
$m$ and contain division points of the factorization of $B_{n,k}$ into
basic blocks $B_{m,j}$.

Let $D(n,k)$ denote the length of $B_{n,k}$ when regarded as a word on
the alphabet whose symbols are basic blocks of level $m$; thus $D(n,k)$ is one more
than the number
of division points in the factorization under discussion. We have
\[
\begin{gathered}
D(m,j)=1 \text{  for all } j=0,1,\dots ,m,\\
D(n,0)=D(n,n)=1 \text{  for all } n \geq m,\\
D(n,k)=D(n-1,k-1)+D(n-1,k) \text{  for all } n>m, k=1, \dots, n-1.
\end{gathered}
\]
Viewing the array whose $(n,k)$ entry is $D(n,k)$ as the sum of $m+1$ Pascal triangles originating at
all the places in level $m$ (since each of these places has an entry
$1$, which is the ``seed'' for an entire Pascal triangle below it), we have
\[
D(n,k)=\sum_{j=0}^mC(n-m,k-j).
\]
By Lemma \ref{lem:quotientbino}, for $n \gg m$ and $k/n \approx \alpha$,
\[
\frac{D(n,k)}{C(n,k)}=\sum_{j=0}^m \frac{C(n-m,k-j)}{C(n,k)} \approx \sum_{j=0}^m
 \alpha ^j(1-\alpha)^{m-j} \leq (\sup \{ \alpha, 1-\alpha \})^{m+1},
\]
so only a small
fraction of the appearances of $B$ in $B_{n,k}$ are ``new''.

The frequency of appearance of $B_{m,j}$ in $B_{n,k}$ is
$C(n-m,k-j)/C(n,k)$,
which has a limiting value as $n \to \infty, k/n \to \alpha$.
Taking into account the number of times that our given block $B$
appears in each basic word $B_{m,j}$ at level $m$, since we can ignore
the ``new'' appearances of $B$ in $B_{n,k}$ we see that $B$ also
has a limiting frequency of appearance in $B_{n,k}$ as $n \to \infty, k/n \to
\alpha$ (and we can compute it, getting the same answer as before).
\end{remark}


\subsection{Topological weak mixing}

In this section we prove that the countable-substitution subshift
version of the Pascal adic is topologically weakly mixing. Since the
system is not minimal (for example, it contains the two fixed points
of $\{ a,b\}^{\mathbb Z}$),
it is not enough to prove that there are no nonconstant continuous
eigenfunctions (cf \cite{petersen3}). Instead, we use a characterization
of topological weak mixing provided by Keynes and Robertson
\cite{keynes-robertson} (see also \cite{keynes}), along with Weyl's theorem on uniform
distribution.

In all the following $(X,T)$ denotes a topological dynamical system,
i.e. $X$ is a compact metric space and $T:X\to X$ is a
homeomorphism. Recall that $(X,T)$ is topologically ergodic
(i.e. topologically transitive) if there is a point $x\in X$ with a
dense orbit.
Since the Pascal adic is ergodic for the Bernoulli
measures $\mu_\alpha$, it follows that the (fully-supported) image measure
$\phi(\mu_\alpha)$ is ergodic for the substitution subshift, and
therefore the system is topologically ergodic. Actually, if $x$
is a path in the Pascal graph which is not eventually diagonal, then
$\phi(x)$ has a dense orbit in $\Sigma$, so in
fact all but countably many orbits in $\Sigma$ are dense.
 We say that $(X,T)$ is {\em topologically
weakly mixing} if $(X \times X, T \times T)$ is topologically ergodic.


\begin{definition} Let $X$ be a topological space. Denote by $\mathcal{C}(X)$ the space of all continuous functions $f:X\to\C$, and by $\B(X)$ the set of bounded functions $f:X\to\C$ such that the set $\mathcal{C}(f)$ of points of continuity of $f$ is residual.
\end{definition}

\begin{definition}
Let $f,g\in\B(X)$. We say that $f$ and $g$ are \emph{essentially equal}, and write $f\ess g$, if $f=g$ on $\mathcal{C}(f)\cap\mathcal{C}(g)$.
\end{definition}

\begin{remark}
If $f,g\in\B(X)$, then $f\ess g$ if and only if $f=g$ on a dense set.
\end{remark}

\begin{thm}[\cite{keynes-robertson}]\label{thm:te}
Let $(X,T)$ be a topological dynamical system. The following are equivalent:
\begin{enumerate}
\item $(X,T)$ is topologically ergodic.
\item For every $f\in\B(X)$, if $f\circ T\ess f$, then $f\ess$ constant.
\item For every $f\in\B(X)$, if $f\circ T=f$ (everywhere), then $f\ess$ constant.
\end{enumerate}
\end{thm}
\begin{remark}\label{rmk:modulusone}
It is easy to show that if $(X,T)$ is ergodic, and $\lambda\in\C$ is an eigenvalue for some $f\in\B(X)$, then $|\lambda|=1$ and $|f|\ess$ constant.
\end{remark}

\begin{definition}
A Borel probability measure $\mu$ on $X$ is called \emph{closed ergodic} for $T$ if every closed invariant subset of $X$ has $\mu$-measure 0 or 1.
\end{definition}



Here is a criterion for weak mixing that we will use to show that the substitution dynamical system $(\Sigma,\sigma)$ is topologically weakly mixing:

\begin{thm}[\cite{keynes-robertson}]\label{thm:twm}
Let $(X,T)$ be a topological dynamical system and suppose there exists a $T$-invariant Borel probability measure $\mu$ supported on all of $X$ which is closed ergodic. The following are equivalent:
\ben
\item $(X,T)$ is topologically weakly mixing.
\item For every $f\in\B(X)$, if there is $\lambda \in \mathbb C$ such that $f\circ T=\lambda f$ (everywhere), then $f\ess$ constant.
\een
\end{thm}




Every measure $\nu_\alpha$ is ergodic (so in particular closed
ergodic) and has support equal to $\Sigma$, so Theorem \ref{thm:twm} applies
in our case. Given $f\in\B(\Sigma)$ and $\lambda \in \mathbb C$ such
that $f\circ \sigma=\lambda f$ (everywhere), if we can show that
$\lambda=1$ then, combining Theorem \ref{thm:te} and Theorem
\ref{thm:twm}, one would show that $(\Sigma,\sigma)$ is topologically
weakly mixing.

We need the following lemma:

\begin{lemma}\label{lem:continuity}
Let $(X,T)$ be a topological dynamical system (with underlying metric $d$). Let $f:X\to\C$ be a function with $\mathcal{C}(f)\ne\emptyset$ and such that $f\circ T=\lambda f$ for some $\lambda\in\C$. Then every point $x\in X$ with dense forward orbit and dense backward orbit is in $\mathcal{C}(f)$.
\end{lemma}

\begin{proof}
Let $x\in X$ be a point such that both semiorbits $\{T^nx|n\ge 0\}$ and $\{T^nx|n\le 0\}$ are dense. Let $z\in\mathcal{C}(f)$. Since $\{T^nx|n\ge 0\}$ is dense, there exists $n_i\nearrow\infty$ such that $T^{n_i}x\to z$.  Therefore $|f\circ T^{n_i}(x)|=|\lambda|^{n_i}|f(x)|$ forces $|\lambda|\le 1$. Since $\{T^nx|n\le 0\}$ is also dense, a similar argument shows that $|\lambda|\ge 1$, so $|\lambda|=1$.

Assume that $x$ has both semiorbits dense and $x\notin\mathcal{C}(f)$.  Then there exists $\alpha>0$ such that for all $\delta>0$ there is $y\in X$ such that
\begin{equation}\label{equ:discont}
d(x,y)<\delta\text{ and } |f(x)-f(y)|\ge\alpha.
\end{equation}
Since $z\in\mathcal{C}(f)$, there exists $\eta>0$ such that
\begin{equation}\label{equ:cont}
\text{if}\quad d(z,u)<\eta, \quad\text{then}\quad |f(z)-f(u)|<\alpha/4.
\end{equation}
Since $x$ has dense orbit we can find $n\in\Z$ such that $d(T^nx,z)<\eta/2$. By continuity of $T^n$ there exists $\delta>0$ such that
\begin{equation*}
\text{if}\quad d(x,y)<\delta, \quad\text{then}\quad |T^nx-T^ny|<\eta/2.
\end{equation*}
For that previous $\delta$, by (\ref{equ:discont}), there exists $y$ such that \[d(x,y)<\delta\text{ and } |f(x)-f(y)|\ge\alpha.\] Since $\lambda$ is an eigenvalue with modulus one we have
\begin{equation*}
|f\circ T^nx-f\circ T^ny|=|\lambda^nf(x)-\lambda^nf(y)|=|f(x)-f(y)|\ge\alpha.
\end{equation*}
On the other hand, since $d(T^nx,z)<\eta/2$ and $d(T^ny,z)<\eta$, by (\ref{equ:cont}) we get the following contradiction: \[|f\circ T^nx-f\circ T^ny|<|f\circ T^nx-f(z)|+|f(z)-f\circ T^ny|<\alpha/2.\]
Thus $x\in\mathcal{C}(f)$.
\end{proof}

\begin{thm}\label{thm:pascal-twm}
The countable substitution dynamical system $(\Sigma,\sigma)$ is topologically weakly mixing.
\end{thm}

\begin{proof}
Let $f\in\B(\Sigma)$ be such that $f\circ \sigma=\lambda f$ for some $\lambda\in \C$. Assume that $\lambda\ne 1$. Since there are no rational eigenvalues
(see \cite{santiago}), $\lambda=e^{2\pi i\beta}$ for some irrational number $\beta$. Suppose that we can find $\omega\in \Sigma$ with the properties that:
\begin{enumerate}
\item $\omega$ has dense orbit,
\item there exist $N_k \to \infty$ such that $\lambda^{N_k}\to -1$,
\item $\sigma^{N_k}\omega\to \omega$.
\end{enumerate}
Then, by Lemma \ref{lem:continuity}, since every point with dense orbit is a continuity point of $f$, the relation $f\circ\sigma^{N_k}\omega=\lambda^{N_k}f(\omega)$ would lead to a contradiction.
To find such a point $\omega$ we use Weyl's theorem on uniform distribution. Since $\beta$ is irrational, for every $k\ge 1$ the set $\{C(n,k)\beta\,:\,n\in\N\}$ is uniformly distributed modulo 1. Therefore there exist $n_1<n_2<\dots<n_k<\dots$ such that \begin{equation}\label{equ:ii}\left|\lambda^{C(n_k,k)}+1\right|<\frac{1}{k} \quad \text{for all }k.\end{equation} Let $x$ be the path in the Pascal graph defined by \[x=01^{n_1}01^{n_2-n_1-1}01^{n_3-n_2-1}0\dots 1^{n_k-n_{k-1}-1}0\dots.\] Since $x$ is not eventually diagonal, $\omega=\phi(x)$ has a dense orbit in $\Sigma$, and hence condition (i) is satisfied. Let $N_k=C(n_k,k)$. Condition (ii) follows from (\ref{equ:ii}), and condition (iii) is guaranteed by the Kink Lemma (Lemma \ref{lem:kink_lemma}).
\end{proof}


\section{Questions}

Many properties of the Pascal adic and related systems remain to be
determined. In particular, the question of weak mixing of the systems
$(X,T,\mu_p)$ remains open. Since ergodicity of the Bernoulli measures
$\mu_p$ under the Pascal adic map $T$ implies the Hewitt-Savage 0,1
Law, and weak mixing is stronger than ergodicity, results along these
lines would have probabilistic implications.

Recently the loose Bernoulli property of each $(X,T,\mu_p)$ has been established by de la Rue and Janvresse \cite{delarue-janvresse}. We believe that these systems have infinite rank, and indeed that they do not have local rank one (see \cite{ferenczi1} for the definitions). Determination of the joinings, or even factors, of these systems would be of considerable interest, as would any more information about their spectra (simple? singular?).

Dynamical properties of a class of adic transformations generalizing the Pascal adic, which code adic transformations on certain shifts of finite type (see \cite{gibbs}), are studied in a forthcoming paper by the first-named author \cite{xman2}.

\begin{ack*}
The first author was partially supported by the Millennium Nucleus in
Information and Randomness, Programa Iniciativa Cientifica Milenio P01-005.
We thank R. Burton, R. Kenyon, and A. Quas for helpful
conversations. This paper is based on the Ph.D. dissertation of the
first author \cite{xman}, written under the direction of the second author.
\end{ack*}


\bibliographystyle{amsplain}



\end{document}